\documentclass[12pt]{amsart}
\usepackage{amsmath,amsfonts,amsthm,amscd,times}
\newtheorem{theorem}{Theorem}
\newtheorem{claim}{Claim}
\newtheorem{proposition}{Proposition}

\newtheorem{definition}{Definition}
\newtheorem{corollary}{Corollary}
\newtheorem{remark}{Remark}
\numberwithin{equation}{section}
\numberwithin{proposition}{section}
\numberwithin{lemma}{section}
\numberwithin{claim}{section}
\numberwithin{corollary}{section}
\newcommand{\bull}{\ensuremath{{}\bullet{}}}

\newcommand{\cpn}{\ensuremath{\mathbb{P}^{N}}}
\newcommand{\slnc}{\ensuremath{SL(N+1,\mathbb{C})}}

\newcommand{\dlb}{\ensuremath{\overline{\partial}}}
\newcommand{\gr}{\ensuremath{\mathbb{G}(N-n-1,N}}
\newcommand{\dl}{\ensuremath{\partial}}
\newcommand{\ra}{\ensuremath{\longrightarrow}}
\newcommand{\ba}{\ensuremath{\begin{align*}}}
\newcommand{\ea}{\ensuremath{\end{align*}}}

\newcommand{\vplt}{\ensuremath{\varphi_{\lambda(t)}}}
\newcommand{\vp}{\ensuremath{\varphi}}

\input{diagrams}

\begin{document}
\title[CM Stability I]{CM Stability and the Generalized Futaki Invariant I}
\author{Sean Timothy Paul}
\email{stpaul@math.wisc.edu}
\address{Mathematics Department at the University of Wisconsin, Madison}
 \thanks{Supported by an NSF DMS grant \# 0505059} 
\author{Gang Tian}            
\email{tian@math.princeton.edu}       
\address{Department of Mathematics, 	 Princeton University}
 \date{April 19, 2008}
\maketitle
 \vspace{-6mm}
\begin{abstract}{Based on the Cayley, Grothendieck, Knudsen Mumford theory of determinants we extend the CM polarization to the Hilbert scheme.   
We identify the weight of this refined line bundle with the generalized Futaki invariant of Donaldson. 
We are able to conclude that CM stability implies K-Stability.
An application of the Grothendieck Riemann Roch Theorem shows that this refined sheaf is isomorphic to the CM polarization introduced by Tian in 1994 on any closed, simply connected base .}
\end{abstract}
\section{Introduction}
\subsection{Statement of Results}  
Throughout this paper $\mathbf{X}$ and $S$ denote complex projective varieties (or schemes) satisfying the following conditions.  
 \smallskip
  \begin{enumerate}
  \item $\mathbf{X}\subset S\times\cpn $; \  $\cpn$  denotes the complex projective space of \emph{lines} in $\mathbb{C}^{N+1}$ .\\
  \item $p_1: \mathbf{X}\rightarrow S$ is flat of relative dimension $n$, degree $d$ with Hilbert polynomial $P$.\\
  \item $L:=p_2^*(\mathcal{O}_{\cpn}(1))$ where $p_2$ is the projection of $\mathbf{X}$ to $\cpn$ .\\
  \item $L|_{\mathbf{X}_z}$ is very ample and the embedding ${\mathbf{X}_z}:=p_1^{-1}(z)\overset{L}{\hookrightarrow} \cpn$ is given by a complete linear system for $z\in S$.
  \end{enumerate}
 \smallskip
 
  It is well known that  $(1)$ and $(4)$ are equivalent to
  \begin{align}
  \mathbb{P}({p_1}_*L)\cong S\times\cpn \ .
  \end{align}
 Which in turn is equivalent to the existence of a line bundle $\mathcal{A}$ on $S$ such that
\begin{align}
{p_1}_*L\cong \underbrace{\bigoplus \mathcal{A}}_{N+1}\ .
 \end{align}

Below $\mbox{{Chow}}(\mathbf{X}\big / S)$ denotes the Chow form of the family $\mathbf{X}\big / S$,
 $\mu$ is the coefficient of $k^{n-1}$ in  $P(k)$, and 
  $\mathcal{M}_{n}$ is the coefficient of $\binom{m}{n}$ in the CGKM expansion of   
 $\det({p_1}_*L^{\otimes m})$ for $m>>0$ . A complete discussion of these notions is given in  section 2.3. 
  We define an invertible sheaf on $S$ as follows. 
  \begin{definition} (The Refined CM polarization \footnote{We use this terminology in order to distinguish this sheaf from one introduced by the first author in (\cite{psc}).}) \\
  \begin{align}\label{defn}
 {\mathbb{L}}_{1}(\mathbf{X}\big / S):= \{\mbox{\emph{Chow}}(\mathbf{X}\big / S)\otimes \mathcal{A}^{d(n+1)}\}^{n(n+1)+\mu}\otimes\mathcal{M}_{n}^{-2(n+1)}
 \end{align}
   \end{definition}

When $S= \mathfrak{Hilb}_{\mathbb{P}^{N}}^{P}(\mathbb{C})$ and $\mathbf{X}= \mathcal{U} \ (\mbox{the universal family})$ we will write $\mathbb{L}_{1}$ instead of $\mathbb{L}_1(\mathbf{X}\big / S)$ .
In our first Theorem we show that the weight of a special linear $\mathbb{C}^{*}$ action on ${\mathbb{L}}^{\vee}_1$ \footnote{$E^{\vee}$ denotes the dual of $E$.}over the Hilbert scheme is the generalized Futaki invariant  of the corresponding degeneration.
The generalized Futaki invariant is defined in section 2.5. In the body of the paper $G$ will always denote $\slnc$ which acts by the standard representation on $\cpn$.
\begin{theorem} (The weight of the Refined CM polarization)\\
\ \\
i) There is a natural $G$ linearization on the line bundle ${{\mathbb{L}}_{1}}$ . \\
\ \\
ii) Let $\lambda$ be a one parameter subgroup of $G$. Let $z\in  \mathfrak{Hilb}_{\mathbb{P}^{N}}^{P}(\mathbb{C})$. Let $w_{\lambda}(z)$ denote the weight  of the restricted $\mathbb{C}^{*}$ action (whose existence is asserted in i)) on 
${{\mathbb{L}}^{\vee}_{1}}|_{z_{0}}$  where  $z_{0}=\lambda(0)z$. Then
\begin{align}\label{weight}
{w_{\lambda}(z)=F_{1}(\lambda)}\ .
 \end{align}
$F_{1}(\lambda)$ is the  generalized Futaki invariant of $\mathbf{X}_{z}$ with respect to $\lambda$.  
\end{theorem}
\begin{remark} We should point out to the reader that given an algebraic group $H$ and an $H$ variety $X$ (or scheme) and an invertible sheaf $L$ on $X$ it is not automatic that the action of $H$ extends to $L$. That is, $L$ need not admit an $H$ linearization.  If $X$ is normal and proper, and $H$ is irreducible then $L$ admits a unique linearization provided that $\chi(H)=\{1\}$ (see \cite{dolga} for a clear account of these facts) . 

\end{remark}
\begin{proposition} 
Assume that the Grothendieck Riemann Roch theorem holds for the map $p_1$. Then the first Chern class of  $c_{1}( {\mathbb{L}}_{1})$ is given by the following formula.
\smallskip
\begin{align}
c_{1}( {\mathbb{L}}_{1})={p_1}_*\left(c_1(K_{\mathbf{X}/ S})c_1({L})^n+\mu\ c_1({L})^{n+1}\right) \quad K_{\mathbf{X} / S}:= K_{\mathbf{X}}\otimes {p^*_1}(K^{\vee}_S)\ .
\end{align}
\end{proposition}
Under the preceding hypothesis the second author introduced the following invertible sheaf on $S$,   the \emph{CM polarization} of the family $\mathbf{X}\overset{p_1}\rightarrow S$. First consider the virtual bundle over $\mathbf{X}$.
\begin{align}
2^n\mathcal{E}_{\mathbf{X} / S}:=(n+1)\left(K^{-1}_{\mathbf{X}}-K_{\mathbf{X}}\right)(L-L^{-1})^n-\mu \left(L-L^{-1}\right)^{n+1} \ .
\end{align}
\begin{definition} (The CM polarization of $\mathbf{X}\overset{p_1}\rightarrow S$)
\begin{align}
\mathbb{L}_{S}:=\mbox{\emph{det}}({p_{1}}_*(\mathcal{E}_{\mathbf{X} / S}))^{\vee} \ .
\end{align}
\end{definition}
\begin{corollary}
Assume that the map 
\begin{align*}
c_1: \mbox{Pic}(S)\rightarrow H^2(S,\mathbb{Z})
\end{align*}
is injective.
Then the refined CM polarization is isomorphic to the CM polarization.
\end{corollary}
\begin{remark} $\mathbb{L}_{S}$ exists  when $\mathbf{X}$ and $S$ are smooth. Whereas ${\mathbb{L}}_{1}(\mathbf{X}\big / S)$ exists for any flat family of schemes. The corollary says that when both exist they are isomorphic, provided the base is, for example, simply connected. 
\end{remark}
\subsection{Resume of results of part II}
In the sequel to this paper we establish the following results. Below we assume that $\mathbf{X}$ and $S$ are smooth projective varieties and that there is an action of $G$ on $\mathbf{X}\big / S$ which restricts to the standard action of $G$ on the fibers $G\ni \sigma:\mathbf{X}_z\rightarrow \sigma \mathbf{X}_z=\mathbf{X}_{\sigma z}$ . $\nu_{\omega}$ denotes the Mabuchi energy.
\begin{theorem} 
 Let $|| \ ||$ be any Hermitian metric on $\mathbb{L}^{\vee}_{1}$. 
 Then there is a continuous function $\Psi_{S}:S\setminus \Delta\rightarrow (-\infty, \ c)$ 
 such that for all $z\in S\big / \Delta$ 
 \begin{align}\label{singnrm}
&d(n+1)\nu_{\omega|_{\mathbf{X}_z}}(\vp_\sigma)= {\log}\left(e^{(n+1)\Psi_{S}(\sigma z)}\frac{||\ ||^{2}(\sigma z)}{||\ ||^{2}( z)}\right)\ .
\end{align}
$c$ denotes a constant which depends only on the choice of background K\"ahler metrics on $S$ and $\mathbf{X}$,  $\Delta$ denotes the discriminant locus of the map $p_1$, and $\omega |_{\mathbf{X}_z}$ denotes the restriction of the Fubini Study form of $\cpn$ to the fiber $\mathbf{X}_z$ .
 \end{theorem}
Under the same hypothesis as the preceding Theorem we have the following corollary of 
(\ref{weight}) and (\ref{singnrm}) .
 \begin{corollary}
 Let $\vplt$ be the Bergman potential associated to an algebraic 1psg $\lambda$ of $G$, and let $z\in S\setminus\Delta$ . Then there is an asymptotic expansion
\begin{align}\label{asymp}
d(n+1){ \nu_{\omega|_{\mathbf{X}_z}}(\varphi_{\lambda(t)})-\Psi_{S}({\lambda(t)})= 
 F_{1}(\lambda)\log(|t|^2)+O(1)} \ \mbox{as}\ |t| \rightarrow 0. \quad 
\end{align}
Moreover  $\Psi_{S}(\lambda(t))= \psi(\lambda)\log(|t|^{2})+O(1)$  where  $\psi({\lambda}) \in \mathbb{Q}_{\geq 0} $ and  $\psi({\lambda}) \in \mathbb{Q}_{+} $ if and only if ${\lambda(0)}\mathbf{X}_z=\mathbf{X}_{\lambda(0)z}$ (the limit cycle  \footnote {See \cite{sopv} pg. 61.} of $\mathbf{X}_z$ under $\lambda$ ) has a component of multiplicity greater than one. $O(1)$ denotes any quantity which is bounded as $|t|\rightarrow 0$.
\end{corollary} 
Moser iteration and a refined Sobolev inequality (see \cite{simon} and \cite{sobolev}) yield the following.
\begin{corollary}
If $\nu_{\omega|_{\mathbf{X}_z}}  $ is proper (bounded from below) then the generalized Futaki invariant of $\mathbf{X}_z$ is {strictly negative} (nonnegative) for all $\lambda \in G$.   \end{corollary}

 \begin{remark} 
We call the left hand side of (\ref{asymp}) the \emph{reduced} K-Energy along $\lambda$.
We also point out that while it is certainly the case that  $F_{1}(\lambda)$ may be defined for any subscheme of $\cpn$ it evidently only controls the behavior of the K-Energy when ${\lambda(0)\mathbf{X}_z} $ is  reduced.
\end{remark}

\begin{remark}
 {The precise constant $d(n+1)$ in front of $\nu_{\omega}$ is not really crucial, since what really matters  is the} sign \emph{of $F_{1}(\lambda)+ \psi({\lambda})$.} 
 {That $\Psi_{S}(\lambda(t))$ has logarithmic singularities can be deduced from \cite{ags}.}
\end{remark}

\begin{remark}
 {We emphasize that we \emph{do not} assume the limit cycle is smooth.}
\end{remark}
 
 \subsection{Background and Motivation from Geometric Analysis}
Let  $(X,\omega)$ be a compact K\"ahler manifold ($\omega$ not necessarily a Hodge class) and $P(X,\omega):= \{\varphi \in  C^{\infty}(X):\omega_{\varphi}:= \omega + \frac{\sqrt{-1}}{2\pi}\dl\dlb\varphi >0\}$  the space of K\"ahler potentials. This is the usual description of all K\"ahler metrics in the same class as $\omega$ (up to translations by constants).  It is not an overstatement to say that the most basic problem in K\"ahler geometry is the following
\begin{center}
\emph{Does there exist  $\varphi \in P(X,\omega)$ such that } $\mbox{Scal}(\omega_{\varphi})\equiv \mu$ . $(*)$\end{center}
This is a fully nonlinear \emph{fourth order} elliptic partial differential equation for $\varphi$. 
  $\mu$ is a constant, the average of the scalar curvature, it depends only on $c_{1}(X)$ and $[\omega]$. When $c_{1}(X)>0$ and $\omega$ represents the \emph{anticanonical} class
 a simple application of the Hodge Theory shows that $(*)$ is equivalent to the \emph{Monge-Ampere equation}.
  \begin{align*}
  \frac{\mbox{det}(g_{i\overline{j}}+\varphi_{i\overline{j}})}{\mbox{det}(g_{i\overline{j}})}=e^{F-\kappa\varphi} \quad ( \kappa=1) \qquad (**)
  \end{align*}
 $F$ denotes the Ricci potential. When $\kappa=0$ this is the celebrated Calabi problem solved by S.T.Yau in the 70's.
  It is well known that $(*)$ is actually a \emph{Variational} problem. There is a natural energy  on the space $P(X,\omega)$ whose critical points are those $\varphi$ such that $\omega   _{\varphi}$ has constant scalar curvature (csc).
This energy was introduced by T. Mabuchi (\cite{mabuchi}) in the 1980's. It is called the \emph{\textbf{K-Energy map}} (denoted by $\nu_{\omega}$) and is given by the following formula
\begin{align*}
 \qquad \nu_{\omega}(\varphi):= -\frac{1}{V}\int_{0}^{1}\int_{X}\dot{\varphi_{t}}(\mbox{Scal}(\varphi_{t})-\mu)\omega_{t}^{n}dt.
\end{align*}
Above, $\varphi_{t}$ is a smooth path in $ P(X,\omega)$ joining $0$ with $\varphi$. The K-Energy does not depend on the path chosen.   
In fact there is the following well known formula for $\nu_{\omega}$ where $O(1)$ denotes a quantity which is bounded on $P(X,\omega)$.
\begin{align*}
&\nu_{\omega}(\varphi)=\int_{X}\mbox{log}\left(\frac{{\omega_{\varphi}}^{n}}{\omega^{n}}\right)\frac{{\omega_{\varphi}}^{n}}{V} - \mu(I_{\omega}(\varphi)-J_{\omega}(\varphi))+O(1) \\
& J_{\omega}(\varphi):= \frac{1}{V}\int_{X}\sum_{i=0}^{n-1}\frac{\sqrt{-1}}{2\pi}\frac{i+1}{n+1}\dl\varphi \wedge \dlb
\varphi\wedge \omega^{i}\wedge {\omega_{\varphi} }^{n-i-1}\\
&I_{\omega}(\varphi):= \frac{1}{V}\int_{X}\varphi(\omega^{n}-{\omega_{\varphi}}^{n})
\end{align*}
It is easy to see that $J_{\omega}(\varphi)\geq 0$ for all $\varphi\in P(X,\omega)$. It is also clear that
\begin{align*}
\frac{1}{n+1}I_{\omega}(\varphi)\leq J_{\omega}(\varphi)\leq \frac{n}{n+1}I_{\omega}(\varphi).
\end{align*}
We have written this down in the case when $\omega = c_{1}(X)$ in which case $\mu=n>0$. In particular observe that $\nu_{\omega}$ is essentially the \emph{difference} of two positive terms. What is of interest for us is that  the problem $(*)$ is not only a {variational} problem but a \emph{minimization} problem. With this said we have the following fundamental result. \ \\

\noindent \textbf{Theorem} (S. Bando and T. Mabuchi \cite{bandmab}) \\
 \emph{If $\omega=c_{1}(X)$ admits a K\"ahler Einstein metric then $\nu_{\omega}\geq 0$. The absolute minimum is taken on the solution to $(**)$ (which is unique up to automorphisms of $X$).}\\
\ \\
Therefore a \emph{necessary} condition for the existence of a K\"ahler Einstien metric is a bound from below on $\nu_{\omega}$. In order to get a \emph{sufficient} condition
one requires that the K-energy \emph{grow} at a certain rate. Precisely, it is required that the K-Energy be \emph{proper}.
This concept was  introduced by the second author  in \cite{psc}.
\begin{definition} 
{$\nu_{\omega}$ is \textbf{proper} if there exists a strictly increasing function $f:\mathbb{R}_{+}\ra\mathbb{R}_{+}$ (where $\lim_{T\ra \infty}f(T)=\infty$) such that $\nu_{\omega}(\varphi)\geq f(J_{\omega}(\varphi))$ for all $\varphi\in P(M,\omega)$}.  
 \end{definition}
 
\noindent \textbf{Theorem} (\cite{psc})\\
  \emph{ Assume that $Aut(X)$ is discrete. Then $\omega=c_{1}(X)$ admits a K\"ahler Einstein metric if and only if $\nu_{\omega}$ is proper}.\\
\ \\
 The next result has recently been established by the second author and Xiuxiong Chen. It holds in an \emph{arbitrary} K\"ahler class $\omega$.\\
\ \\
 \textbf{Theorem} (\cite{chentian})\\
  \emph{ If $\omega$ admits a metric of csc then} $\nu_{\omega}\geq 0$.\\
\subsection{Geometric Invariant Theory and K-Energy Asymptotics}  
 The motivation for our (purely algebraic) work is to analyze the behavior of $\nu_{\omega}$  along a large but finite dimensional group $G$  of \emph{matrices} in the polarized case. That is, we assume that $\omega=c_1(L)$ where $L$ is a very ample line bundle over $X$ which provides an embedding of $X$ into $\cpn$. It is well known that there is a map 
 from $G$ into $P_{\omega}(X)$ given by $G\ni\rightarrow \log\frac{|\sigma z|^2}{|z|^2} \in P_{\omega}(X)$. Therefore, the K-Energy map may be considered as a map from $G$ to $\mathbb{R}$.
 It is when we restrict $\nu_{\omega}$ to $G$ that we make the connection with Mumfords' Geometric Invariant Theory. The past couple of years have witnessed quite a bit of activity in K\"ahler geometry due to this connection. For example, the reader is invited to consult , as a rather small sample, the following papers \cite{skdtoric,skdproj,geodrays,ross&thomas1, ross&thomas2,gauduchon,mabuchiext,mabuchichow,harmonicchern}.   
\section{The Refined CM Polarization}
\subsection{Hilbert Points and the Numerical Criterion}
The purpose of this section is to define precisely the weight of certain $\mathbb{C}^{*}$ actions on the Hilbert point of an \emph{individual} projective variety. Later we extend this notion to \emph{families}. 

Let $(X,L)$ be a polarized algebraic variety. Assume that $L$ is very ample with associated embedding
 \begin{align*}
 X\underset{\varphi_{L}}\longrightarrow \mathbb{P}(H^{0}(X,L)^{*}).
 \end{align*}
 Fix  an isomorphism
 \begin{align*}
 \sigma:H^{0}(X,L)^{*}\underset{\cong}\longrightarrow \mathbb{C}^{N+1}\ .
 \end{align*}
 In this way we consider $X$ embedded in $\cpn$.
 Let $m\in \mathbb{Z}$ be a large positive integer.
 Then there is a surjection\footnote{$\mathbf{S}^{m}$ denotes the $m$th symmetric power operator.}
 \begin{align*}
 \Psi_{X,m}:\mathbf{S}^{m}({\mathbb{C}^{N+1}})^{\vee}\longrightarrow H^{0}(X,\mathcal{O}(m))\rightarrow 0\ .
 \end{align*}
 Let $P(m)=P(X,\mathcal{O}(m))=h^{0}(X,\mathcal{O}(m))$ and $d_m:= \mbox{dim}(\mathbf{S}^{m}({\mathbb{C}^{N+1}}))$ . It is a deep fact (see \cite{eckart}) that there is an integer $m(P)$ depending \emph{only} on the Hilbert polynomial $P$ such that for all $m\geq m(P)$, the kernel of $ \Psi_{X,m}$ 
 \begin{align*}
 \mbox{Ker}(\Psi_{X,m})\in \mbox{G}(P(m),\mathbf{S}^{m}({\mathbb{C}^{N+1}})^{\vee}) 
 \end{align*}\footnote{The Grassmannian of $P(m)$ dimensional \emph{quotients} of $\mathbf{S}^{m}({\mathbb{C}^{N+1}})^{\vee}$.}
 completely determines $X$. In other words,  the \emph{entire} homogeneous ideal can be recovered from its $m$th graded piece for all subschemes $X$ of $\cpn$ with Hilbert polynomial $P$ .
 We have the Pl\"ucker embedding
\begin{align*}
 \mathcal{P}:\mbox{G}(P(m),\mathbf{S}^{m}({\mathbb{C}^{N+1}})^{*}) \rightarrow \mathbb{P}\left(\bigwedge^{d_{m}-P(m)}\mathbf{S}^{m}(\mathbb{C}^{N+1})^{*}\right).
 \end{align*}
Next we consider the canonical nonsingular  pairing 
 \begin{align*}
 \bigwedge^{d_{m}-P(m)}\mathbf{S}^{m}(\mathbb{C}^{N+1})^{*}\otimes \bigwedge ^{P(m)}\mathbf{S}^{m}(\mathbb{C}^{N+1})^{*}\longrightarrow \mbox{\bf{det}}(\mathbf{S}^{m}(\mathbb{C}^{N+1})^{\vee})\ .
 \end{align*}
 This induces a natural isomorphism
 \begin{align*}
 \mathbb{P}\left(\bigwedge^{d_{m}-P(m)}\mathbf{S}^{m}(\mathbb{C}^{N+1})^{\vee}\right)\overset{\iota}\cong \mathbb{P}\left(\bigwedge^{P(m)}\mathbf{S}^{m}(\mathbb{C}^{N+1})\right).
 \end{align*}
 Combining this identification with the Pl\"ucker embedding we associate to $\mbox{Ker}(\Psi_{X,m})$ a unique point, called (following Gieseker) the $m$th Hilbert Point
 \begin{align*}
\mbox{Hilb}_{m}(X):= \iota(\mathcal{P}(\mbox{Ker}(\Psi_{X,m})))\in  \mathbb{P}\left(\bigwedge^{P(m)}\mathbf{S}^{m}(\mathbb{C}^{N+1})\right).
\end{align*}
 We can give a coordinate description of this as follows.
  Given the surjection
 \begin{align*}
 \Psi_{X,m}:\mathbf{S}^{m}({\mathbb{C}^{N+1}})^{\vee}\longrightarrow H^{0}(X,\mathcal{O}(m))\rightarrow 0
 \end{align*}
 we can take its \emph{determinant}
 \begin{align*}
 \wedge^{P(m)}\Psi_{X,m}:\wedge^{P(m)} \mathbf{S}^{m}({\mathbb{C}^{N+1}})^{\vee}\rightarrow \mathbf{det}H^{0}(X,\mathcal{O}(m))\cong \mathbb{C}.
 \end{align*}
 We associate to this map a point in the projective space \emph{dual to} $\wedge^{P(m)} \mathbf{S}^{m}({\mathbb{C}^{N+1}})$
 \begin{align*}
 [\wedge^{P(m)}\Psi_{X,m}]\in \mathbb{P}(\wedge^{P(m)} \mathbf{S}^{m}({\mathbb{C}^{N+1}})^{\vee}).
 \end{align*}
 We then have that
 \begin{align*}
 [\wedge^{P(m)}\Psi_{X,m}]=\mbox{Hilb}_{m}(X).
 \end{align*}
 Let $\mathbf{I}=(i_{0},i_{1},\dots,i_{N})$ be a multiindex with $|\mathbf{I}|:=i_{0}+i_{1}+\dots +i_{N}=m$, $i_{j}\in \mathbb{N}$. Let $e_{0},e_{1},\dots, e_{N}$ be the standard basis of $\mathbb{C}^{N+1}$, and $z_{0},z_{1},\dots,z_{N}$ be the dual basis of linear forms. Consider the monomials
 $M_{\mathbf{I}}:= e_{0}^{i_{0}}e_{1}^{i_{1}}\dots e_{N}^{i_{N}}$ and $M^{*}_{\mathbf{I}}:= z_{0}^{i_{0}}z_{1}^{i_{1}}\dots z_{N}^{i_{N}}$. Fix a basis $\{f_{1},\dots, f_{P(m)}\}$ of $H^{0}(X,\mathcal{O}(m))$. 
 Then 
\begin{align*}
\begin{split}
\wedge^{P(m)}\Psi_{X,m}(M^{*}_{{\mathbf{I}}_{j_{1}}}\wedge \dots \wedge M^{*}_{{\mathbf{I}}_{j_{P(m)}}})=&\Psi_{X,m}(j_{1},\dots , j_{P(m)})f_{1}\wedge \dots \wedge  f_{P(m)} \\
 &\Psi_{X,m}(j_{1},\dots , j_{P(m)})\in \mathbb{C}\ .
\end{split}
\end{align*}
Then in homogeneous coordinates we can write
\begin{align*}
\wedge^{P(m)}\Psi_{X,m}=\sum_{\{{\mathbf{I}}_{j_{1}},\dots,{\mathbf{I}}_{j_{P(m)}}\}}\Psi_{X,m}(j_{1},\dots , j_{P(m)})    M_{{\mathbf{I}}_{j_{1}}}\wedge \dots \wedge M_{{\mathbf{I}}_{j_{P(m)}}}.
 \end{align*}
 \subsection{Weights}
Let $\lambda:\mathbb{C}^{*}\rightarrow G$ be an algebraic one parameter subgroup. We may assume that $\lambda$ has been diagonalized on the standard basis $\{e_{0},e_{1},\dots, e_{N}\}$. Explicitly, we assume that there are  $r_{i}\in \mathbb{Z}$ such that
 \begin{align*}
 \lambda(t)e_{j}=t^{r_{j}}e_{j} \quad \mbox{with the obviously necessary property}\quad \sum_{0\leq j\leq N}r_{j}=0.
 \end{align*}
Define the weight of $\lambda$ on the monomial $M_{\mathbf{I}}$ by\footnote{Of course this makes sense even if $\lambda$ is not special linear.}
\begin{align*}
w_{\lambda}(M_{\mathbf{I}}):= r_{0}i_{0}+r_{1}i_{1}+\dots +r_{N}i_{N}.
\end{align*}
Now we make the following
\begin{definition} (Gieseker \cite{globmod})
The  \emph{\textbf{weight}} of the $m$th Hilbert point of $X$ is the integer 
\begin{align*}
w_{\lambda}(m):=\overset{\mbox{Min}}{\{{\mathbf{I}}_{j_{1}},\dots,{\mathbf{I}}_{j_{P(m)}}\}} \left(\sum_{1\leq k\leq P(m)}w_{\lambda}(M_{{\mathbf{I}}_{j_{k}}})|\Psi_{X,m}(j_{1},\dots , j_{P(m)})\neq 0\right)\ .
\end{align*}
\end{definition}
Concerning this $w_{\lambda}(m)$ Mumford proved the following (non obvious) result in \cite{sopv} :\\
\ \\
 \noindent \emph{For fixed $\lambda$ and large $m$ the weight of the Hilbert point is given by a numerical polynomial of degree at most $n+1$}
\begin{align*}
w_{\lambda}(m)=a_{n+1}(\lambda)m^{n+1}+a_{n}(\lambda)m^{n}+O(m^{n-1}).
\end{align*}

We may consider, by slight abuse, $\mbox{Hilb}_{m}(X):= \wedge^{P(m)}\Psi_{X,m}$ as a point in a \emph{vector space}, namely
\begin{align*}
\mbox{Hilb}_{m}(X)\in \wedge^{P(m)} \mathbf{S}^{m}({\mathbb{C}^{N+1}})\big / \{0\}.
\end{align*}
A moments thought shows that the  weight  of the $\mathbb{C}^{*}$ action $\lambda$  is  the unique integer $w_{\lambda}(m)$ such that
\begin{align*}
{\lim_{t\rightarrow 0}t^{-w_{\lambda}(m)}\lambda(t)\mbox{Hilb}_{m}(X) \quad \mbox{\emph{exists and is not equal to zero}}}.
\end{align*}
 Sometimes we will write $w_{\lambda}( \mbox{Hilb}_{m}(X))$ in place of $w_{\lambda}(m)$ to indicate the dependence on the underlying  complex projective manifold (or scheme) $X$.

Next, endow $\wedge^{P(m)} \mathbf{S}^{m}({\mathbb{C}^{N+1}})$ with \emph{any} hermitian metric $||\  ||$, for example the standard one. Then a reformulation of what we have just said runs as follows. We may characterize the weight as the leading term of the small $t$ asymptotics of the \emph{logarithm} of the norm
 \begin{align*}
 {\log(||\lambda(t)\mbox{Hilb}_{m}(X)||^{2})=w_{\lambda}(m)\log(t^2)+O(1)}\ .
 \end{align*}
 
 This characterisation will be relevant in the sequel to this paper .  

\begin{definition} (D. Gieseker \cite{globmod})\\
$\mbox{Hilb}_{m}(X)$ is (\textbf{semi}) {\textbf{stable}} provided $w_{\lambda}(m)(\leq 0)<0$ for all one parameter subgroups $\lambda$. 
\end{definition}
We must mention that in his \emph{tour de force} work (\cite{globmod}) D. Gieseker verified the stability of Hilbert points of  pluricanonical models of algebraic surfaces of general type.
 
\subsection{Hilbert Points and Stability in Families}
 The preceding notions may be extended to the relative setting.  
 As stated in the introduction $\mathbf{X}$ and $S$ denote complex projective varieties (or schemes) satisfying the following conditions.  
 \smallskip
  \begin{enumerate}
  \item $\mathbf{X}\subset S\times\cpn $; \  $\cpn$  denotes the complex projective space of \emph{lines} in $\mathbb{C}^{N+1}$ .\\
  \item $p_1: \mathbf{X}\rightarrow S$ is flat of relative dimension $n$, degree $d$ and Hilbert polynomial $P$.\\
  \item $L:=p_2^*(\mathcal{O}_{\cpn}(1))$ where $p_2$ is the projection of $\mathbf{X}$ to $\cpn$ .\\
   \item $L|_{\mathbf{X}_z}$ is very ample and the embedding ${\mathbf{X}_z}:=p_1^{-1}(z)\overset{L}{\hookrightarrow} \cpn$ is given by a complete linear system for $z\in S$.
  \end{enumerate}
 \smallskip
  $(1)$ is equivalent to
  \begin{align}\label{projtriv}
  \mathbb{P}({p_1}_*L)\cong S\times\cpn \ .
  \end{align}
   (\ref{projtriv}) is equivalent to the existence of a line bundle $\mathcal{A}$ on $S$ such that
\begin{align}
{p_1}_*L\cong \underbrace{\bigoplus \mathcal{A}}_{N+1}\ .
 \end{align}
Let $m$ be an integer large enough so that  $R^{i}{p_1}_{*}(L^{\otimes m})=0$  when $ i>0$. In this case ${p_1}_{*}(L^{\otimes m})$ is a locally free sheaf  of rank $r_{m}:=P(m)$ ,  and we may define the determinant
\begin{align*}
\mbox{det}({p_1}_{*}L^{\otimes m}):=\bigwedge^{r_{m}}({p_1}_{*}L^{\otimes m}).
\end{align*}  
Observe that
\begin{align*}
\mbox{det}({p_1}_{*}L)\cong \mathcal{A}^{P(1)}  \ . 
\end{align*}
For $m\gg 0$ we have a surjective map 
\begin{align*}
\bigwedge^{P(m)}\mathbf{S}^{m}({p_1}_{*}L)\rightarrow \mbox{det}({p_1}_{*}L^{m})\rightarrow 0\ .
\end{align*}
Which in turn yields the map 
\begin{align*}
\mathcal{A}^{mP(m)}\otimes \bigwedge^{P(m)}\mathbf{S}^{m}
(\underset{N+1}{\underbrace{\bigoplus}} \mathcal{O}_{S})\rightarrow \mbox{det}({p_1}_{*}L^{m})\rightarrow 0 \ .
\end{align*}
Throwing $\mathcal{A}^{mP(m)}$ onto the other side gives
\begin{align*}
\bigwedge^{P(m)}\mathbf{S}^{m}
(\underset{N+1}{\underbrace{\bigoplus}} \mathcal{O}_{S})\rightarrow \mbox{det}({p_1}_{*}L^{m})\otimes\mathcal{A}^{-mP(m)}\rightarrow 0 \ .
\end{align*}
This latter is equivalent to
\begin{align*}
\bigwedge^{P(m)}\mathbf{S}^{m}
(\underset{N+1}{\underbrace{\bigoplus}} \mathcal{O}_{S})\rightarrow \mbox{det}({p_1}_{*}L^{m})\otimes \mbox{det}({p_1}_{*}L)^{\frac{-mP(m)}{N+1}}\rightarrow 0 \ .
\end{align*}
From which we deduce the existence of a map $\varphi_{m}$ into the Hilbert scheme.
\begin{align*}
\begin{split}
&\varphi_{m}:S\rightarrow \mathfrak{Hilb}_{\mathbb{P}^{N}}^{P}(\mathbb{C})\hookrightarrow \mathbb{P}^{N(m)}:= \mathbb{P}(\bigwedge^{P(m)}\mathbf{S}^{m}
(\underset{N+1}{\underbrace{\bigoplus}} \mathbb{C})) \\ 
&\varphi_{m}^{*}\mathcal{O}_{\mathbb{P}^{N(m)}}(1)\cong \mbox{det}({p_1}_{*}L^{m})\otimes \mbox{det}({p_1}_{*}L)^{\frac{-mP(m)}{N+1}} \ .
\end{split}
\end{align*}
Then the \emph{Hilbert polarization of the family} $\mathbf{X}\overset{p_1}\rightarrow S$ is by definition, the following invertible sheaf on $S$
\begin{align*}
{\mbox{Hilb}_{m}(\mathbf{X}/S) :=  \mbox{det}({p_1}_{*}L^{m})\otimes \mbox{det}({p_1}_{*}L)^{\frac{-mP(m)}{N+1}}}\cong \mbox{det}({p_1}_{*}L^{m})\otimes \mathcal{A}^{-mP(m)}\ .
\end{align*}
  We note that this sheaf $\mbox{Hilb}_{m}(\mathbf{X}/S)$ is $G$ linearized
 in the case when the family comes equipped with a $G$ action satisfying the requirements listed in the introduction. Let $z\in S$, since $S$ is closed there is a point $z_{0}\in S$ such that $\lambda(0)z:= lim_{t\rightarrow 0}\lambda(t)z=z_{0}$. This gives a \emph{one dimensional} representation of $\mathbb{C}^{*}$
 \begin{align*}
  \mbox{Hilb}_{m}(\mathbf{X}/S)^{(-1)}|_{z_{0}}.
\end{align*}
So that  $\mathbb{C}^{*}$ acts by $t^{w_{z}(m)},\ w_{z}(m) \in \mathbb{Z}$. Then it is easy to see that

\begin{align*}
w_{z}(m)= w_{\lambda}( \mbox{Hilb}_{m}(X_{z}))=w_{\lambda}( \mbox{Hilb}_{m}(X_{z_{0}})). 
\end{align*}
\subsection{Numerical Functions}

By a \emph{numerical function} we mean simply any mapping 
\begin{align*}
\chi:\mathbb{Z}\rightarrow \mathbb{Z} \ .
\end{align*}
Let $\chi \in \mathbb{Q}[T]$ have degree $k$, then we have the well known fact (see J.P. Serre \cite{ localalg})
\begin{align*}
&\chi(m)\in \mathbb{Z} \ \mbox{ \emph{for all $m$ if and only if }} 
\chi(m)=\sum_{0\leq i \leq k}a_{i}\binom{m}{i} \ \mbox{\emph{where}} \  a_{i}\in \mathbb{Z}.
\end{align*}

Let $\chi:\mathbb{Z}\rightarrow \mathbb{Z}$ be a numerical function recall that $\chi$ is \emph{eventually polynomial} provided that the following holds:\\
\ \\
\emph{There is a polynomial $P\in \mathbb{Q}[T]$ and an integer $m_{0}$ such that for all $m\geq m_{0}$ we have }
\begin{align*}
\chi(m)=P(m).
\end{align*}
For any numerical function $\chi$ we define the standard \emph{forward difference operator} $\Delta$ by the formula
\begin{align*}
\Delta \chi (m):= \chi(m+1)-\chi(m)\ .
\end{align*}
By induction we have
\begin{align*}
\Delta^{k+1}\chi(m)=\sum_{0\leq i \leq k+1} (-1)^{i+1}\binom{k+1}{i}\chi(m+i) \ .
\end{align*}
 
\begin{proposition}(see \cite{localalg})\label{serre}
Let  $\chi:\mathbb{Z}\rightarrow \mathbb{Z}$ be a numerical function. Then the following are equivalent.
\begin{align*}
& i) \ \chi \ \mbox{is eventually polynomial}.\\
\ \\
&ii)\  \mbox{There is an integer $k\in \mathbb{Z}_{+}$ such that $\Delta^{k} \chi (m)\equiv e_{k}(\chi)$ (a constant)}.\\
\ \\
&iii)\  \mbox{There is an integer $k\in \mathbb{Z}_{+}$ such that } \\ &\chi(m)=\sum_{0\leq i \leq k}a_{i}\binom{m}{i} \ \mbox{where}\  a_{i}\in \mathbb{Z};
\mbox{and we have that}\ \Delta^{k} \chi (m) = a_{k}.
\end{align*}
\end{proposition}

Following Grothendieck, Knudsen, and Mumford we apply the previous considerations on numerical functions not only  to $\chi(m)$ but to the\emph{ Picard group valued} numerical function
\begin{align*}
m \in \mathbb{Z}_{+}\rightarrow \mbox{det}({p_1}_{*}L^{m}) \in \mbox{Pic}(S).
\end{align*}

We denote by $\Delta^{k+1} \mbox{det}({p_1}_{*}L^{\otimes m})$ the \emph{difference sheaves} defined inductively
\begin{align*}
&\Delta \mbox{det}( {p_1}_*L^{m}):= \mbox{det}( {p_1}_{*}L^{ m+1})\otimes \mbox{det}({p_1}_*L^{ m})^{(-1)}\\
\ \\ 
&\Delta^{2}\mbox{det}( {p_1}_{*}L^{m}):= \Delta^{}\mbox{det}({p_1}_{*}L^{ m+1})\otimes  \Delta^{}\mbox{det}({p_1}_*L^{ m})^{(-1)}\\
\ \\
&\ \qquad \qquad \qquad \qquad \dots \\
&\Delta^{k+1} \mbox{det}({p_1}_{*}L^{ m})\cong \bigotimes_{i=0}^{k+1}\mbox{det}({p_1}_{*}L^{ m+i})^{(-1)^{i+1}\binom{k+1}{i}}.
\end{align*}
Let $\mathfrak{C}(n,d; \cpn)$ denote the \emph{Chow Variety}  of dimension $n$  and degree $d$ algebraic cycles inside $\cpn$.  $\mathfrak{C}(n,d; \cpn)$ is known to be a subvariety
\footnote{It is by definition a \emph{subset} of this projective space, that it is an actual sub\emph{variety} is a fundamental Theorem of Chow and van der Waerden. }
 of  the projective space of sections of degree $d$ on the Grassmannian.
\begin{align*}
\mathfrak{C}(n,d; \cpn) \overset{\iota}\hookrightarrow \mathbb{P}(H^{0}(\mathbb{G},\mathcal{O}(d))):=\mathbb{P}(H^{0}(\mathbb{G}(N-n-1,\cpn),\mathcal{O}(d)))\  .
\end{align*}

There is a map \footnote{The construction of this map is carried out in the appendix.} (see \cite{git} and \cite{fog})  $\Delta$ from $ \mathfrak{Hilb}_{\mathbb{P}^{N}}^{P}(\mathbb{C})$ to the Chow variety  which sends a subscheme $\mathcal{I} $ of $\cpn$ with Hilbert polynomial $P$ to the Chow form of the top dimensional component of its underlying cycle. For $m$ sufficiently large, we let  $\varphi_{m}$  denote the map from $S$ to $ \mathfrak{Hilb}_{\mathbb{P}^{N}}^{P}(\mathbb{C})$  . Combining this map with $\Delta$ yields the sequence
\begin{align*}
S\overset{\varphi_{m}}\longrightarrow \mathfrak{Hilb}_{\mathbb{P}^{N}}^{P}(\mathbb{C})\overset{\Delta}\longrightarrow \mathfrak{C}(n,d; \cpn) \overset{\iota}\hookrightarrow 
\mathbb{P}(H^{0}(\mathbb{G},\mathcal{O}(d))).
\end{align*}
Let $\mathcal{O}(1)$ denote the hyperplane line on $\mathbb{P}(H^{0}(\mathbb{G},\mathcal{O}(d)))$.
Then we define the \emph{Chow form of the map} $\mathbf{X}\overset{p_1}\rightarrow S$ to be the invertible sheaf on $S$
\begin{align*}
\mbox{Chow}(\mathbf{X} /S ):= \varphi_{m}^{*}{\Delta}^{*}\iota^{*}\mathcal{O}(1).
\end{align*}

In their paper (\cite{detdiv}, Theorem 4)  Knudsen and Mumford \footnote{Following suggestions in an unpublished letter of Grothendieck to Mumford (1962). This result was already known to Cayley in special cases in 1850.} proved the following fundamental results (In fact, much more is true . The family $\mathbf{X}/S$ need \emph{not} be flat, we remark on the more general construction in the appendix . ) \\
 
\noindent I) \emph{On the base $S$ there is a canonical isomorphism of invertible sheaves} 
 \begin{align}\label{knmd}
 \mbox{Chow}(\mathbf{X}/S)\otimes \mathcal{A}^{d(n+1)} \cong \Delta^{n+1}\mbox{det}({p_1}_{*}L^{m})\ .
 \end{align}

\begin{remark} The isomorphism not only says that the right hand side of (\ref{knmd}) is constant, but that  this constant is the Chow form of the map $p_1$.
\end{remark}
In view of Proposition \ref{serre}
we may expand $\mbox{det}({p_1}_{*}L^{m})$ in terms of the binomial polynomials:\\
\ \\
\noindent II) \emph{There are invertible sheaves $\mathcal{M}_{0},\mathcal{M}_{1},\dots, \mathcal{M}_{n+1} $ on $S$ and a canonical and functorial isomorphism}:
 \begin{align*}
 &\mbox{det}({p_1}_{*}L^{m})\cong \bigotimes _{j=0}^{n+1}\mathcal{M}_{j}^{\binom{m}{j}}\ .
  \end{align*}
\emph{ Moreover, the leading term is related to the Chow form of the map $\mathbf{X}\overset{p_1}\rightarrow S$ {as follows}}

\begin{align*}  
  &   \mathcal{M}_{n+1} \cong \mbox{Chow}(\mathbf{X}/S)\otimes \mathcal{A}^{d(n+1)}\ .
 \end{align*}
 Now we restate our definition.\\
  
 {\textbf{Definition}} (\emph{The Refined CM polarization}) \\
 \begin{align}\label{defn}
 {\mathbb{L}}_{1}(\mathbf{X}\big / S):= \{\mbox{\emph{Chow}}(\mathbf{X} / S)\otimes \mathcal{A}^{d(n+1)}\}^{n(n+1)+\mu}\otimes\mathcal{M}_{n}^{-2(n+1)}
 \end{align}
\subsection{The weight of the Refined CM Polarization}
In this subsection we show that the weight of ${\mathbb{L}}_{1}(\mathbf{X}\big / S)$ with respect to $\lambda$ is the generalized Futaki invariant $F_1(\lambda)$ .To carry this out it is very convenient to construct an entire \emph{sequence} of sheaves ${\mathbb{L}}_{l}(\mathbf{X}\big / S)$ for $0\leq l \leq n+1$. These sheaves are then identified with certain polynomial combinations of the sheaves $\mbox{Hilb}_{m+i}(\mathbf{X}/S)$. The weights of the sheaves ${\mathbb{L}}_{l}(\mathbf{X}\big / S)$ are then computed by appealing to the expansion of the weight of $\mbox{Hilb}_{m+i}(\mathbf{X}/S)$. Then we let $l=1$ to complete the proof of our main theorem.

Fix any subvariety $X\ra \cpn$ with Hilbert polynomial $P$.
We begin by considering the ratio.
 \begin{align*}
&\frac{w_{\lambda}(m)}{mP(m)} =F_{0}(\lambda)+F_{1}(\lambda)\frac{1}{m}+\dots +F_{l}(\lambda)\frac{1}{m^{l}}+\dots\\
&\mbox{Observe that we may expand the coefficient of $m^{-l}$ as follows}\\
&F_{l}(\lambda)= c_{l,n+1}a_{n+1}(\lambda)+c_{l,n}a_{n}(\lambda)+c_{l,n-1}a_{n-1}(\lambda)+\dots + c_{l,n+1-l}a_{n+1-l}(\lambda)\\
&\mbox {where the $c_{l,j}$ are all rational functions of the coefficients of the Hilbert polynomial $P$}\ . 
\end{align*}

\begin{definition}(Donaldson (\cite{skdtoric}))\newline
$F_{1}(\lambda)$ is the generalized Futaki invariant of $X$ with respect to $\lambda$ .
\end{definition}
Now we may introduce invertible sheaves $\mathbb{L}_{l}(\mathbf{X}/S) $ on $S$ for all $l=0,1,2,\dots,n+1$ as follows.
 
\begin{definition}  
 \begin{align}
\mathbb{L}_{l}(\mathbf{X}/S):=\bigotimes_{k=n+1-l}^{n+1}{\mathcal{M}_{k}}^{\frac{1}{k!}\sum_{0\leq j\leq k-1}(-1)^{j+1}\sigma_{j}(1,2,\dots,k-1)c_{l.k-j}}
\end{align}
Where the $\mathcal{M}_{k} ,\ 0\leq k \leq n+1$ are the coefficients in the Cayley, Grothendieck, Knudsen, Mumford expansion.
 \end{definition}
Let $f_{l}(m):=m^{l}$. Then we define polynomials $P_{k,l}(m)$
\begin{align*}
P_{k,l}(m):=\Delta^{k}f_{l}(m) \ . 
\end{align*}
It is easy to see that
\begin{align*}
P_{k,l}(m)=\sum_{0\leq j \leq k}(-1)^{j+1}\binom{k}{j}(m+j)^{l}\ .
\end{align*}
It is not difficult to verify that
\begin{align*}
P_{k,l}(m)=
\begin{cases}
(-1)^{k+1}k! , & \mbox{if} \ k=l \\
0,& \mbox{if} \  l<k \ .
\end{cases}
\end{align*}
In general, $P_{k,k+d}(m)$ is a polynomial in $m$ of degree $d$.
Given $0\leq l \leq n+1$ let 
\begin{align*}
(q_{n+1}(m), q_{n}(m), q_{n-1}(m), \dots, q_{n+1-l}(m))
\end{align*}
be the unique solution to the equation
\begin{align*}
\begin{pmatrix} P_{n+1,n+1}(m)&P_{n,n+1}(m)&\dots & \dots & P_{n+1-l,n+1}(m)\\
      0&P_{n,n}(m)& P_{n-1,n}(m)&\dots & P_{n+1-l,n}(m)  \\
       0& 0&P_{n-1,n-1}(m)  &\dots& P_{n+1-l,n-1}(m)\\
       0&  0& 0& \dots &\dots  \\ 
       \dots&\dots&\dots&\dots\\
         0&  0& 0& \dots&  P_{n+1-l,n+1-l}(m)
       \end{pmatrix}
\begin{pmatrix} q_{n+1}(m)\\
q_{n}(m)\\
q_{n-1}(m)\\
\dots \\
\dots\\
q_{n+1-l}(m)
\end{pmatrix}
=
 \begin{pmatrix} c_{l,n+1}\\
c_{l,n}\\
c_{l,n-1}\\
\dots \\
\dots\\
c_{l,n+1-l}
\end{pmatrix}
\end{align*}
Then we have the following proposition.
\begin{proposition} For all $l=0,1,2 \dots, n+1$ we have
\begin{align}\label{sum0}
\bigotimes_{0\leq p\leq l}\bigotimes_{0\leq i \leq n+1-p} \mbox{\emph{Hilb}}_{m}(\mathbf{X} / S)^{(-1)^{i}q_{n+1-p}(m)\binom{n+1-p}{i}} \cong  \mathbb{L}_{l}(\mathbf{X}/S)
\end{align}
 \end{proposition} 
 \begin{corollary}If the family $\mathbf{X}/S$ admits an action of $\slnc$ then
 the sheaves $\mathbb{L}_{l}(\mathbf{X}/S)$ admit natural linearizations.
 \end{corollary}
 \begin{proof}
Writing out the left hand side of (\ref{sum0}) gives 
\begin{align*}
&\bigotimes_{0\leq p\leq l}\bigotimes_{0\leq i \leq n+1-p}\mathbf{det}(f_{*}L^{m+i})^{(-1)^{i}q_{n+1-p}(m)\binom{n+1-p}{i}}\otimes \mathcal{A}^{\frac{1}{N+1}(-1)^{i+1}(m+i)P(m+i)q_{n+1-p}(m)\binom{n+1-p}{i}}\ . \\
\end{align*}
The exponent of $\mathcal{A}$  satisfies the following
\begin{align}\label{sum}
\sum_{0\leq p\leq l}\sum_{0\leq i\leq n+1-p}(-1)^{i+1}(m+i)P(m+i)q_{n+1-p}(m)\binom{n+1-p}{i}=\begin{cases}1,&l=0\\
0,& l>0\ .
\end{cases}
\end{align}
 
 To see this we first write
 \begin{align*}
 (m+i)P(m+i)=b_{n}(m+i)^{n+1}+b_{n-1}(m+i)^{n}+\dots + b_{j}(m+i)^{j+1}+\dots
\end{align*}
Then the left hand side of (\ref{sum}) is given by
\begin{align}\label{sum2}
\sum_{0\leq p\leq l}\sum_{0\leq i\leq n+1-p}\sum_{0\leq j \leq n}(-1)^{i+1}  q_{n+1-p}(m)b_{j}(m+i)^{j+1}\binom{n+1-p}{i}
\end{align}
Recall that we have defined the polynomials $P_{n+1-p,j+1}(m)$ by the formula
\begin{align}\label{P}
P_{n+1-p,\ j+1}(m)= \sum_{0\leq i \leq n+1-p}(-1)^{i+1}(m+i)^{j+1}\binom{n+1-p}{i}\ .
\end{align}
Substituting (\ref{P}) into (\ref{sum2}), switching the order of summation  and appealing to the definiton of the $q_{k}(m)$ gives
\begin{align}\label{sum3}
\sum_{n-l\leq j \leq n}\sum_{n-j\leq p\leq l}q_{n+1-p}(m)P_{n+1-p,\ j+1}(m)b_{j}=\sum_{n-l\leq j \leq n}b_{j}c_{l,\ j+1}.
\end{align}
By definition of the $c_{l,k}$ the right hand side of (\ref{sum3}) is the coefficient of ${m}^{-l}$ in the expansion of
\begin{align*}
\frac{mP(m)}{mP(m)}\equiv 1\ .
\end{align*}
From now on we will assume that $l>0$. With this assumption we have
\begin{align}\label{sum4}
\begin{split}
 &\bigotimes_{0\leq p\leq l}\bigotimes_{0\leq i \leq n+1-p} \mbox{{Hilb}}_{m+i}(\mathbf{X}/ S)^{(-1)^{i}q_{n+1-p}(m)\binom{n+1-p}{i}} \cong \\
&\bigotimes_{0\leq p\leq l}\bigotimes_{0\leq i \leq n+1-p}\mathbf{det}({p_1}_{*}L^{m+i})^{(-1)^{i}q_{n+1-p}(m)\binom{n+1-p}{i}} \cong \\
&  \bigotimes_{0\leq p\leq l}\bigotimes_{0\leq i \leq n+1-p}\bigotimes_{0\leq k \leq n+1}\mathcal{M}_{k}^{(-1)^{i}q_{n+1-p}(m)\binom{n+1-p}{i}\binom{m+i}{k}}\cong\\
&\bigotimes_{0\leq k \leq n+1}\mathcal{M}_{k}^{\sum_{0\leq p\leq l}\sum_{0\leq i \leq n+1-p }(-1)^{i}q_{n+1-p}(m)\binom{n+1-p}{i}\binom{m+i}{k}}.
\end{split}
\end{align}
Next we study the exponent of $\mathcal{M}_{k}$ on the last line of (\ref{sum4}).
First we expand the binomial coefficients in of powers of $(m+i)$ (below $\sigma_j$ denotes the $jth$ elementary symmetric function ) .
 \begin{align*}
 \binom{m+i}{k}=\frac{1}{k!}\sum_{j=0}^{k-1}(-1)^{j}\sigma_{j}(1,2,\dots,k-1)(m+i)^{k-j}.
 \end{align*}
 So that we have
 \begin{align*}
&\sum_{0\leq i \leq n+1-p }(-1)^{i}\binom{n+1-p}{i}\binom{m+i}{k}=\\
&\sum_{0\leq j\leq k-1}\frac{1}{k!}(-1)^{j}\sigma_{j}(1,2,\dots,k-1) \sum_{0\leq i \leq n+1-p}(-1)^{i}\binom{n+1-p}{i}(m+i)^{k-j}=\\
&\sum_{0\leq j\leq k-1}\frac{1}{k!}(-1)^{j+1}\sigma_{j}(1,2,\dots,k-1) P_{n+1-p,k-j}(m)\ .
 \end{align*}
Therefore,
\begin{align*}
&\sum_{0\leq p\leq l}\sum_{0\leq i \leq n+1-p }(-1)^{i}q_{n+1-p}(m)\binom{n+1-p}{i}\binom{m+i}{k}=\\
&\sum_{0\leq p\leq l}\sum_{0\leq j\leq k-1}\frac{1}{k!}(-1)^{j+1}\sigma_{j}(1,2,\dots,k-1) P_{n+1-p,k-j}(m)q_{n+1-p}(m)=\\
&\frac{1}{k!}\sum_{0\leq j\leq k-1}(-1)^{j}\sigma_{j+1}(1,2,\dots,k-1)c_{l.k-j}.
\end{align*}
Which completes the proof of the proposition.  
\end{proof}
Below we assume that $\mathbf{X}$ and $S$ are complex varieties (or schemes) and that there is an action of $G$ on $\mathbf{X}\big / S$ which restricts to the standard action of $G$ on the fibers $G\ni \sigma:\mathbf{X}_z\rightarrow \sigma \mathbf{X}_z=\mathbf{X}_{\sigma z}$ .
 In this situation we have the following identity of weights.
\begin{align}\label{wts}
w_{\lambda}(\mathbb{L}^{\vee}_{l}(\mathbf{X}/S))=w_{\lambda}\left(\bigotimes_{0\leq p\leq l}\bigotimes_{0\leq i \leq n+1-p} \mbox{{Hilb}}_{m+i}(\mathbf{X}/ S)^{(-1)^{i+1}q_{n+1-p}(m)\binom{n+1-p}{i}}\right)\ .
\end{align}
 \begin{claim}
 The weight on the right hand side of (\ref{wts}) is given by $F_{l}(\lambda)$. 
\end{claim}
\begin{proof}
Since the weight is additive under tensor product the right hand side is given by the sum
\begin{align*}
&\sum_{0\leq j\leq n+1}\sum_{0\leq p\leq l}\sum_{0\leq i \leq n+1-p }(-1)^{i+1}q_{n+1-p}(m)a_j(\lambda)(m+i)^j\binom{n+1-p}{i} \\
&= \sum_{0\leq j \leq n+1}\sum_{0\leq p\leq l}P_{n+1-p, j}(m)q_{n+1-p}(m)a_{j}(\lambda) \\
&= \sum_{n+1-l \leq j \leq n+1}c_{l,j}a_{j}(\lambda) \\
&=F_{l}(\lambda)\ .
\end{align*}
Now let $l=1$.
\end{proof}
\subsection{The first Chern Class of $\mathbb{L}_1(\mathbf{X}/S)$}\ \newline
In this section we assume that the Grothendieck Riemann Roch theorem holds for the map $\mathbf{X} \overset{p_1}\rightarrow S$ .
On $\mathbf{X}$ we introduce the virtual bundles for $0\leq l \leq n+1$.
\begin{align*}
\mathcal{E}_{l}(m):= \sum_{\{0\leq j \leq l\}}\sum_{\{0\leq i \leq n+1-j\}}q_{n+1-j}(m)(-1)^{i}\binom{n+1-j}{i}L^{m+i}  \ .
\end{align*}
Then it is easy to see that
\begin{align*}
 &\bigotimes_{0\leq p\leq l}\bigotimes_{0\leq i \leq n+1-p} \mbox{{Hilb}}_{m+i}(\mathbf{X}/ S)^{(-1)^{i}q_{n+1-p}(m)\binom{n+1-p}{i}} \cong \mathbf{det}({p_1}_*\mathcal{E}_{l}(m))\ .
 \end{align*}
Therefore we have 
\begin{align}
\mathbb{L}_1(\mathbf{X}/S)\cong \mathbf{det}({p_1}_*\mathcal{E}_{1}(m))\ .
\end{align}
Therefore,
\begin{align}
c_{1}(\mathbb{L}_1(\mathbf{X}/S)) = c_1(\mathbf{det}({p_1}_*\mathcal{E}_{1}(m)))\ .
\end{align} 
The computation of the Chern character  of $\mathcal{E}_{l}(m)$ follows the same pattern as the computations in the preceding sections and is left to the reader. This completes the proof. $\Box$
\subsection{Appendix on Det and Div: The Cayley, Grothendieck, Knudsen, Mumford Expansion}\footnote{There is a closely related work of Fogarty (\cite{fog}) on this subject .} \newline
In this section we outline the main application of the determinant construction from \cite{detdiv} to Cayley-Chow forms.
For complete proofs the reader should consult this paper.
Historically, the construction is due to Cayley (see \cite{cay}) who showed that any \emph{resultant} can be expressed as the determinant of a \emph{complex}. Grothendieck extended Cayley's idea to show that the Cayley-Chow form of \emph{any} projective variety may be expressed as the determinant of a complex. As we will recall in what follows Grothendieck emphasized that one really takes the Chow form of  \emph{coherent sheaf} on $\cpn$ of dimension $n$ (see \cite{eisenbud} ) .

Let $X$ be a complex variety (or scheme). 
\begin{definition}
A complex $(E^{\bull},\ \dl_{\bull})$ of (coherent) sheaves on $X$  is \textbf{perfect} provided it is locally quasiisomorhpic to a bounded complex of vector bundles.
\end{definition}
 Let $(E^{\bull},\ \dl_{\bull})$ and $F^{\bull}$ be two perfect complexes on $X$. Assume that there is a proper subvariety $W$ of $X$ and a quasiisomorphism $\iota$
\begin{align*}
(E^{\bull},\ \dl_{\bull})\overset{{\iota}}{\cong}  (F^{\bull}, \ \delta_{\bull}) \ \mbox{over}\ X\setminus W \ .
\end{align*}
When $(F^{\bull}, \ \delta_{\bull})=0^{\bull}$ we say that $(E^{\bull},\ \dl_{\bull})$ is \emph{generically exact} .
In this situation Grothendieck, Knudsen, and Mumford associate to this data a \emph{Cartier} divisor $Div(\iota)$ on $X$ such that
\begin{align}
\mathcal{O}(Div(\iota))\cong Det(E^{\bull})^{\vee}\otimes Det(F^{\bull})
\end{align}
In the case when $(E^{\bull},\ \dl_{\bull})$ is generically exact we write $Div(E^{\bull})$ in place of $Div(0^{\bull}\rightarrow E^{\bull})$.

The following is a reformulation of "condition $Q_k$ " on page 50 from ( \cite{detdiv} ) into more classical language.
\begin{definition}
Let $\mathbf{X}\overset{p}\rightarrow S$ be a proper map between complex varieties. Let $(E^{\bull},\ \dl_{\bull})$ be a perfect complex on $\mathbf{X}$. Then we say that $(E^{\bull},\ \dl_{\bull})$ satisfies the condition $Q_k$ relative to the map $p$ provided the following conditions are met .\\
\ \\
i) There is a subvariety $W\subset S$ {such that for all} $s\in S\setminus W$
\begin{align*}
\mbox{dim}(Supp(E^{\bull},\ \dl_{\bull})\cap \mathbf{X}_{s})\leq k 
\end{align*}
ii) For all hypersurfaces $Z\subset S$ we have
\begin{align*}
\mbox{dim}(Supp(E^{\bull},\ \dl_{\bull})\cap \mathbf{X}_{s})\leq k +1 \ \mbox{generically on $Z$} \ .
\end{align*}
\end{definition}
Recall that $Supp(E^{\bull},\ \dl_{\bull})$ denotes the locus of points in $\mathbf{X}$ where the complex fails to be exact.

 Let  $\mathbf{X}\overset{p}\rightarrow S$ be a proper map of finite Tor-dimension between complex varieties .\footnote{Recall that $p$ has \emph{finite Tor-dimension} provided there exists a finite resolution of $p_*\mathcal{O}_{\mathbf{X}}$ by flat $\mathcal{O}_{S}$ modules.} Grothendieck has shown (see Proposition 4.8 SGA 6, expose 3 LNM 225, p. 257) that $R^{\bull}p_*(E^{\bull},\ \dl_{\bull})$ is perfect on $S$ whenever $(E^{\bull},\ \dl_{\bull})$ is perfect on $\mathbf{X}$ .
\begin{proposition}
Let  $\mathbf{X}\overset{p}\rightarrow S$ be a proper map of finite Tor-dimension between complex varieties .
Let $(E^{\bull},\ \dl_{\bull})$ be a perfect complex satisfying condition $Q_{-1}$ for the map $p$. Then\\
\ \\
i) $R^{\bull}p_*(E^{\bull},\ \dl_{\bull})$ is generically exact on $S$ and we may define $Div(R^{\bull}p_*(E^{\bull},\ \dl_{\bull}))$.\\
\ \\
ii) For all line bundles $L$ on $\mathbf{X}$ we have
\begin{align*}
Div(R^{\bull}p_*(E^{\bull},\ \dl_{\bull}))=Div(R^{\bull}p_*(E^{\bull}\otimes L,\ \dl_{\bull})) \ .
\end{align*}
\end{proposition}
The principal application of these ideas is to construct Cayley-Chow forms in families. An immediate by product of which is the existence of the refined CM polarization $\mathbb{L}_1(\mathbf{X}/S)$ for \emph{non flat} families.

Let $\mathcal{E}$ be a vector bundle of rank $N+1$ on $S$. Let $\mathbb{P}=\mathbb{P}(\mathcal{E})$ denote the corresponding projective bundle, and $\widehat{\mathbb{P}}=\mathbb{P}(\mathcal{E}^{\vee})$ the dual projective bundle. We form the fiber square.
\begin{diagram}
 \mathbb{P}\times_{S}\widehat{\mathbb{P}}^{(n+1)} &\rTo^{p_2}&\widehat{\mathbb{P}}^{(n+1)}\\
\dTo ^{p_{1}}&&\dTo^{\widehat{\pi}}\\
\mathbb{P} &\rTo ^{\pi}& S\\
 \end{diagram}
 There is a canonical section 
 \begin{align}
 \delta \in \Gamma(\mathbb{P}\times_{S}\widehat{\mathbb{P}}^{(n+1)}, \overbrace{\bigoplus \mathcal{O}_{\mathbb{P}}(1)\otimes\mathcal{O}_{\widehat{\mathbb{P}}}(1)}^{n+1})
\end{align}
Let $Z$ denote the zero locus of $\delta$. If we denote the direct sum of $\mathcal{O}_{\mathbb{P}}(1)\otimes\mathcal{O}_{\widehat{\mathbb{P}}}(1)$ by $Q$ then we have the Koszul complex $(K^{\bull}_{(n+1)},\  \delta)$  of sheaves on 
$\mathbb{P}\times_{S}\widehat{\mathbb{P}}^{(n+1)}$ 
\begin{align}
0\rightarrow \bigwedge^{n+1}Q^{\vee}\rightarrow \bigwedge^nQ^{\vee}\rightarrow \dots \rightarrow Q^{\vee}\rightarrow \mathcal{O}_{\mathbb{P}\times_{S}\widehat{\mathbb{P}}^{(n+1)}}\rightarrow \mathcal{O}_{Z}\rightarrow 0
\end{align}
In the next proposition $E^{\bull}(m)$ denotes the complex twisted $m$ times by $\mathcal{O}_{\mathbb{P}}(1)$.
\begin{proposition}
Let $(E^{\bull},\ \dl_{\bull})$ be a perfect complex on $\mathbb{P}$ satisfying condition $Q_n$ for the map $\pi$.
Then $L^{\bull}p^*_1(E^{\bull}(m),\ \dl_{\bull})\overset{L}{\otimes} (K^{\bull}_{(n+1)},\  \delta)$ satisfies condition $Q_{-1}$ for the map $p_2$.
\end{proposition}
 Therefore we can define the \emph{Chow divisor of the complex} $(E^{\bull},\ \dl_{\bull})$ .
 \begin{definition}
 $Chow((E^{\bull},\ \dl_{\bull}))=Div(R^{\bull}{p_2}_*\left(L^{\bull}p^*_1(E^{\bull}(m),\ \dl_{\bull})\overset{L}{\otimes} (K^{\bull}_{(n+1)},\  \delta)\right))$ .
 \end{definition}
 This is a divisor on $\widehat{\mathbb{P}}^{(n+1)}$  satisfying
 \begin{align*}
 \mathcal{O}(Chow((E^{\bull},\ \dl_{\bull})))=\mathcal{O}(Chow((E^{\bull}(m),\ \dl_{\bull})))=det(R^{\bull}{p_2}_*\left(L^{\bull}p^*_1(E^{\bull}(m),\ \dl_{\bull})\overset{L}{\otimes} (K^{\bull}_{(n+1)},\  \delta)\right)) \ .
 \end{align*}
An interesting computation shows that
\begin{align}
det(R^{\bull}{p_2}_*\left(L^{\bull}p^*_1(E^{\bull}(m),\ \dl_{\bull})\overset{L}{\otimes} (K^{\bull}_{(n+1)},\  \delta)\right))\cong  \widehat{\pi}^*\Delta^{n+1}det(R^{\bull}{\pi_*}(E^{\bull}(m),\ \dl_{\bull}))\otimes \mathcal{H}^{\Delta^n\chi(E^{\bull}(m))}
\end{align}
In the formula above $\mathcal{H}:=\overbrace{\bigotimes \mathcal{O}_{\widehat{\mathbb{P}}}(-1)}^{n+1}$.
Since the left hand side is independent of $m$, so is the right hand side. Therefore we have the following main result of 
( \cite{detdiv} ) . \newline
\begin{center}\emph{There are invertible sheaves $\mathcal{M}_j$ $0\leq j \leq n+1$ on $S$ such that}
\begin{align*}
det(R^{\bull}{\pi_*}(E^{\bull}(m),\ \dl_{\bull})\cong \otimes_{j=0}^{n+1}\mathcal{M}_j^{\binom{m}{j}}\  .
\end{align*}
  \emph{Moreover $\chi(E^{\bull}(m))$ is a polynomial of degree at most $n$} .  \end{center}
 
\subsection{Indication of proof when $S$ is a point }
The references for this subsection are (\cite{cay}, \cite{fischer}, \cite{detdiv}, \cite{fog}, \cite{ksz}, \cite{gkz}, and \cite{weyman}). The underlying idea behind the Chow form is to describe
an arbitrary variety by a \emph{single} equation.
Let $X \subset \cpn$ be an $n$ dimensional irreducible subvariety of $\cpn$ with degree $d$, then the Chow form,
or associated hypersurface to $X$ is defined by
\begin{align*}\label{chow}
Z_{X}:= \{L \in \mathbb{G}:= \mathbb{G}(N-n-1,\mathbb{C}P^{N}): L\cap X \neq \emptyset\}.
\end{align*}
It is easy to see that $Z_{X}$ is an \emph{irreducible} hypersurface (of degree $d$) in $\mathbb{G}$.
Since the homogeneous coordinate ring of the grassmannian is a UFD, any codimension one subvariety with
degree $d$ is given by the vanishing of a section  $R_{X}$ of the homogeneous coordinate ring\footnote{See \cite{tableaux}  pg. 140 exercise 7.}
 \begin{align*}
 \{\ R_{X}=0\ \}= Z_{X}  \ ; \  R_{X} \in \mathbf{P}H^{0}(\mathbb{G},\mathcal{O}(d)).
 \end{align*}

By abuse of terminology, we will often call $R_{X}$ the Chow form (or Chow point) of $X$.
Following \cite{ksz} we can be more concrete as follows.
Let $M_{n+1, N+1}^{0}(\mathbb{C})$ be the (Zariski open and dense) subspace
of the vector space of $(n+1)\times(N+1)$ matrices consisting of matrices of full rank.
We have the canonical projection
\begin{align*}
p:M_{n+1, N+1}^{0}(\mathbb{C})\rightarrow \gr ,
\end{align*}
defined  by taking the kernel of the linear transformation. This map is dominant, so the closure of the  preimage
\begin{align*}
\overline{p^{-1}(Z_{X})}\subset \overline{M_{n+1, N+1}^{0}(\mathbb{C})}=
M_{n+1, N+1}(\mathbb{C})
\end{align*}
is also an irreducible hypersurface of degree $d$ in $M_{n+1, N+1}(\mathbb{C})$.
Therefore, there is a unique\footnote{Unique up to scaling.} (symmetric multihomogeneous) polynomial (which will also be denoted by $R_{X}$)
such that 
\begin{align*}
Z:= \overline{p^{-1}(Z_{X})}= \{R_{X}(w_{ij})=0\} \ ;  \ R_{X}(w_{ij})\in \mathcal{P}^{d}[M_{n+1, N+1}(\mathbb{C})].
 \end{align*}
We view the entries $w_{ij}$ as the coefficients of the $n+1$ linear forms
$l_{i}$ defining the plane in $\gr$

\begin{align*}
\begin{pmatrix}l_{0}\\
               \dots \\
               l_{n}\\ \end{pmatrix}
=
\begin{pmatrix}w_{00} & \dots & w_{0N}\\
        w_{10} & \dots & w_{1N}  \\
         \dots& \dots & \dots \\
         w_{n0} & \dots & w_{nN}   \\ \end{pmatrix}\ .
\end{align*}

The purpose of what follows is to provide an (in principal) \emph{explicit formula} for the polynomial $R_{X}(w_{ij})$, which is essentially due to Cayley in his remarkable 1848 note \cite{cay} on resultants. The first rigorous proof seems to be \cite{fischer}. 
We are in debt to the basic paper \cite{detdiv}, as we have said repeatedly throughout this note.  

To begin let $\left(E^{\bull},\dl_{\bull}\right)$ be a bounded complex of finite dimensional $\mathbb{C}$ vector spaces
\[ \begin{CD} 
0@>>>E^{0}@>\dl_{0}>>E^{1}@>\dl_{1}>> \dots @>>> E^{i}@>\dl_{i}>> E^{i+1}@>>> \dots @>\dl_{k-1}>> E^{k}@>>>0 \ . \end{CD}
\]
Recall that the determinant $\mathbf{Det}(E^{\bull})$  of the complex  $\left(E^{\bull},\dl_{\bull}\right)$ is defined to be the one dimensional vector space
\begin{align*}
  \mathbf{Det}(E^{\bull}):=  \bigotimes_{i=0}^{k}\mathbf{det}(E^{i})^{(-1)^{i+1}}\ .
\end{align*}
As usual, for any vector space $V$ we set $V^{-1}:= \mbox{Hom}_{\mathbb{C}}(V,\mathbb{C})$, the dual space to
$V$.
Let $H^{i}(E^{\bull},\dl_{\bull})$ denote the $i^{th}$ cohomology group of this complex.
When $V = \mathbf{0}$ , the zero vector space , we set $\mathbf{det}(V):= \mathbb{C}$.
The determinant of the cohomology is defined in exactly the same way:

\begin{align*}
\mathbf{Det}(H^{\bull}(E^{\bull},\dl_{\bull})):= \bigotimes_{i=0}^{k}\mathbf{det}(H^{i}(E^{\bull},\dl_{\bull}))^{(-1)^{i+1}}\ .
\end{align*}

We have the fundamental facts (\cite{detdiv}):\\
{\bf{D1}}
{Assume that} $(E^{\bull},\dl_{\bull})$ and $(F^{\bull},\delta_{\bull})$ are \emph{quasi-isomorphic}. Then
\begin{align*}
\mathbf{Det}(E^{\bull})\cong \mathbf{Det}(F^{\bull})\ .
\end{align*}
Important corollaries of this fact are the following.\\
{\bf{D2}}
 \emph{There is a canonical isomorphism between the determinant of the complex and}
 \indent \emph{ the determinant of its cohomology}:
\begin{align*}
\tau(\dl_{\bull}):\mathbf{Det}(E^{\bull}) \cong \mathbf{Det}(H^{\bull}(E^{\bull},\dl_{\bull}))\ .
\end{align*}
{\bf{D3}}
\emph{Assume that the complex} $\left(E^{\bull},\dl_{\bull}\right)$ \emph{is acyclic, then} $\mathbf{Det}(E^{\bull})$ \emph{is canonically trivial}:

\begin{align*}
\tau(\dl_{\bull}):\mathbf{Det}(E^{\bull}) \cong \underline{\mathbb{C}}\ .
\end{align*}
It is ({\bf{D3}}) which is relevant for our purpose. It says is that 
\emph{there is a canonically given \textbf{nonzero} element of} $\mathbf{Det}(E^{\bull})$, provided this complex is exact.
The essential ingredient in the formation of the Chow point consists in identifying this canonical ``section''. 
 In order to proceed, we recall the \emph{Torsion} (denoted by $\mathbf{Tor}\left(E^{\bull},\dl_{\bull}\right)$) of the complex $\left(E^{\bull},\dl_{\bull}\right)$. Define $n_{i}:= \mbox{dim}(\dl(E_{i}))$, now choose $S_{i}\in \wedge^{n_{i}}(E_{i})$ with $\dl{S_{i}}\neq 0$, then $\dl{S_{i}}\wedge S_{i+1}$ spans ${\bf{det}}(E_{i+1})$ ( since the complex is exact), that is

\begin{align*}
{\bf{det}}(E_{i+1})= \mathbb{C}\dl{S_{i}}\wedge S_{i+1} .
\end{align*}
With this said we define
\begin{align*}
{{\mathbf{Tor}\left(E^{\bull},\dl_{\bull}\right):= (S_{0})^{-1}\otimes(\dl S_{0}\wedge S_{1})\otimes (\dl S_{1}\wedge S_{2})^{-1}\otimes \dots \otimes (\dl S_{k-1})^{(-1)^{k+1}}}} \ .
\end{align*}
Then we have the following reformulation of {\bf{D3}}.\\ 
\noindent{\bf{D4}}
\begin{align*}
 \mathbf{Tor}\left(E^{\bull},\dl_{\bull}\right) \mbox{\emph{is independent of the choices}}\ S_{i}.
 \end{align*}

By fixing a basis $\{f_{1,i},f_{2,i}, \dots f_{d_{i},i}\}$ in each of the $E_{i}$ ($\mbox{dim}(E_{i}):= d_{i}$), 
we may associate to this \emph{based} exact complex a \emph{scalar}:

\begin{align*}
\mathfrak{Tor}\left(E^{\bull},\dl_{\bull}; \{f_{1,i},f_{2,i}, \dots f_{d_{i},i}\}\right)\in \mathbb{C}.
\end{align*}
 Which is defined through the identity:
\begin{align*}
\mathbf{Tor}\left(E^{\bull},\dl_{\bull}\right)= \mathfrak{Tor}\left(E^{\bull},\dl_{\bull}; \{f_{1,i},f_{2,i}, \dots f_{d_{i},i}\}\right)\mathbf{det}(\dots f_{1,i},f_{2,i}, \dots f_{d_{i},i}\dots).
\end{align*}
Where we have set
\begin{align*}
\mathbf{det}(\dots f_{1,i},f_{2,i}, \dots f_{d_{i},i}\dots):= (f_{1,0}\wedge \dots \wedge 
f_{d_{0},0})^{-1}\otimes \dots \otimes  (f_{1,k}\wedge \dots \wedge f_{d_{k},k})^{(-1)^{k+1}}.
\end{align*}
When we have fixed a basis of our (exact) complex (that is, a basis of each term in the complex) we will call $\mathfrak{Tor}\left(E^{\bull},\dl_{\bull}; \{f_{1,i},f_{2,i}, \dots f_{d_{i},i}\}\right)$ the Torsion of the (based exact) complex.

With these preliminaries out of the way, we can return to the setting of projective geometry. So let  $X$ be an irreducible 
$n$ dimensional subvariety of $\cpn$. 
Let $\mathcal{O}(1)$ be the positive generator of $\mbox{Pic}(\cpn)$. Let $\mathcal{O}(m)_{X}$ be the sheaf of the mth twist of $\mathcal{O}(1)$ restricted to $X$, 
and let $H^{0}(X,\mathcal{O}(m))$ be the space of global sections.
Consider the \emph{Koszul Complex} of sheaves on $X$:
\begin{align*}
(\mathfrak{K}^{\bull}(m)_{X},\dl_{\bull})\qquad \dl_{i}: \mathfrak{K}^{i}(m)_{X}\rightarrow \mathfrak{K}^{i+1}(m)_{X}.
\end{align*}
Where the \emph{ith} term is defined by
\begin{align*}
 \mathfrak{K}^{i}(m)_{X}:= \mathcal{O}(m+i)_{X}\otimes {\wedge^{i} \mathbb{C}^{n+1}}\ .
\end{align*}
The boundary operator is given by
\begin{align*}
\dl(f\otimes \alpha):= \sum_{j=0}^{n}l_{j}f\otimes e_{j}\wedge \alpha \ .
\end{align*}

The $e_{j}$ are the standard basis of $\mathbb{C}^{n+1}$ and the $\{l_{0},\dots,l_{n}\}$ are $n+1$ linear forms on $\cpn$ (\cite{ksz}). In particular we consider a complex whose boundary operators  \emph{depend on parameters}.
We have the following 
\begin{proposition}
The complex $(\mathfrak{K}^{\bull}(m)_{X},\dl_{\bull})$ is exact provided the $n+1$ linear forms 
$\{l_{0},\dots, l_{n}\}$ have no common root on $X$.
\end{proposition}

Now assume that ${m}\gg 0$, then Serres' Theorems \cite{faisceaux} tell us that the higher cohomology groups of the $\mathfrak{K}^{\bull}(m)_{X}$ are all zero, and if the Koszul complex  of \emph{sheaves} is exact then the complex of \emph{global sections} is also exact (\cite{bredon}). In other words we have an acyclic complex of \emph{finite dimensional} vector spaces
\begin{align}\label{complex}
(\Gamma(\mathfrak{K}^{\bull}(m)_{X}), \dl_{\bull}) = (H^{0}(X,\mathcal{O}(m+\bull))\otimes \wedge^{\bull} \mathbb{C}^{n+1},\dl_{\bull}).
\end{align} 

 We can apply the theory \textbf{D1}-\textbf{D4} to such complexes.
So in this case we can, by choosing bases $\{f_{1,i},f_{2,i}, \dots f_{d_{i},i}\}$ of the $\Gamma(\mathfrak{K}^{\bull}(m)_{X}) $
introduce the \emph{Torsion of the Koszul complex}:

\begin{align*}
\mathfrak{Tor}(\Gamma(\mathfrak{K}^{\bull}(m)_{X}),\dl_{\bull};\{f_{1,i},f_{2,i}, \dots f_{d_{i},i}\})(l_{0},\dots, l_{n}).
\end{align*}

This is by construction a \emph{\textbf{rational}} function\footnote{A priori, $\mathfrak{Tor}$ is an alternating product of determinants of certain minors of the matrices representing the boundary operators wrt the bases $\{f_{1,i},f_{2,i}, \dots f_{d_{i},i}\}$.} on $M_{n+1, N+1}(\mathbb{C})$, which is well defined and \emph{\textbf{nowhere zero}} away from $Z= \{R_{X}=0\}$.

Geometrically we have a complex of vector bundles $(\mathcal{E}^{\bull},\ \delta_{\bull})$ over the affine space $\mathbb{V}:=M_{n+1, N+1}(\mathbb{C})$ which is \emph{exact off of a codimension one subvariety}  $Z$. Such a complex is said to be \emph{generically exact}. Therefore there is a nonvanishing \emph{section}  $\sigma$
\begin{align*}
\mathbb{V}\big / Z \overset{\sigma}\rightarrow Det(\mathcal{E}^{\bull},\delta_{\bull}).
\end{align*}
Then
\begin{align*}
\sigma(w_{ij})= \mathfrak{Tor}(\Gamma(\mathfrak{K}^{\bull}(m)_{X}),\dl_{\bull};\{f_{1,i},f_{2,i}, \dots f_{d_{i},i}\})(w_{ij}).
\end{align*}

Since $Z$ is irreducible $\mathfrak{Tor}$ is forced to have either zeros or poles along $Z$. In other words one knows \emph{apriori} that there is an integer $p$ (the $Z$-adic order) such that:

\begin{align*}
\mathfrak{Tor}(\Gamma(\mathfrak{K}^{\bull}(m)_{X}),\dl_{\bull};\{f_{1,i},f_{2,i}, \dots f_{d_{i},i}\})(l_{0},\dots, l_{n}) = R_{X}((l_{0},\dots, l_{n}))^{p}  
\end{align*}
The sensitive point is to determine this order. To begin we \emph{augment} the complex (\ref{complex}) to the following
free complex over the polynomial ring $\mathbb{C}[w_{ij}]$:
\begin{align*}
(H^{0}(X,\mathcal{O}(m+\bull))\otimes{\wedge^{\bull} \mathbb{C}^{n+1}}\otimes \mathbb{C}[w_{ij}],\dl_{\bull}).
\end{align*}
Then we can assert the following
\begin{proposition} The localisation 
\begin{align*}
(H^{0}(X,\mathcal{O}(m+\bull))\otimes{\wedge^{\bull} \mathbb{C}^{n+1}}\otimes \mathbb{C}(w_{ij}),\dl_{\bull})
\end{align*}
of the augmented complex is exact.
\end{proposition}
Algebraically, we are in the following situation. We have a bounded finite free complex $(\mathcal{F}^{\bull},\delta_{\bull})$ over a Noetherian unique factorisation domain $R$
such that its localisation (at zero)  $(\mathcal{F}^{\bull}\otimes_{R}R_{(0)},\delta_{\bull})$  is exact. In the terminology of commutative algebra such complexes are said to be generically exact. Choosing bases $\{f_{1,i},f_{2,i}, \dots f_{d_{i},i}\}$ \emph{over} $R$
we compute the torsion, just as before. We have 
\begin{align*}
\mathfrak{Tor}(\mathcal{F}^{\bull},\delta_{\bull};\{f_{1,i},f_{2,i}, \dots f_{d_{i},i}\})\in {R^{*}}_{(0)}\ \mbox{modulo}\ \mathbb{U}.
\end{align*}
Where $\mathbb{U}$ is the unit group of the ring $R$.
Since $R$ is a UFD we may decompose the determinant into its unique factorisation into irreducibles
\begin{align*}
\mathfrak{Tor}(\mathcal{F}^{\bull},\delta_{\bull};\{f_{1,i},f_{2,i}, \dots f_{d_{i},i}\})=\prod_{f \in Irred(R)}f^{ord_{f} ({det}(\mathcal{F}^{\bull}\otimes_{R}R_{(0)},\delta_{\bull}))}.
\end{align*}
Of course, $ord_{f} ({Det}(\mathcal{F}^{\bull}\otimes_{R}R_{(0)},\delta_{\bull}))=0$ for all but finitely many $f$.
If $f$ is an irreducible element in $R$ then, as usual, $R_{(f)}$ denotes the localisation of $R$ at the prime $f$. Then we have the following (\cite{detdiv} Theorem 3 part $vi)$)\ .
\begin{proposition}\label{order}\ \\
\begin{align*}
&i) \ \mbox{The homology modules} \ \mathbb{H}^{i}(\mathcal{F}^{\bull}\otimes_{R}R_{(f)},\delta_{\bull})\ \mbox{have finite length over}\  R_{(f)}.\\
\ \\
&ii) \ ord_{f} (\mathbf{Det}(\mathcal{F}^{\bull}\otimes_{R}R_{(0)},\delta_{\bull}))=\sum_{i \geq 0}(-1)^{i}l_{R_{(f)}}(\mathbb{H}^{i}(\mathcal{F}^{\bull}\otimes_{R}R_{(f)},\delta_{\bull})).\\
\end{align*}
 Where $l_{R_{(f)}}(\mathbb{H}^{i}(\mathcal{F}^{\bull}\otimes_{R}R_{(f)},\delta_{\bull}))$ denotes the length of the module.
\end{proposition}
\emph{Proof}\\
$i)$ By the universal coefficient theorem we have
\begin{align*}
\mathbf{0}\cong \mathbb{H}^{i}(\mathcal{F}^{\bull}\otimes_{R}R_{(0)},\delta_{\bull}))\cong \mathbb{H}^{i}(\mathcal{F}^{\bull}\otimes_{R}R_{(f)},\delta_{\bull}))\otimes_{R_{(f)}}R_{(0)}\oplus \mbox{Tor}_{1}^{R_{(f)}}(\mathbb{H}^{i-1},R_{(0)})\ .
\end{align*}
Since $R_{(0)}$ is flat over $R_{(f)}$ we have that  $\mbox{Tor}_{1}^{R_{(f)}}(\mathbb{H}^{i-1},R_{(0)})=\mathbf{0}$. Therefore 
\begin{align*}
 \mathbb{H}^{i}(\mathcal{F}^{\bull}\otimes_{R}R_{(f)},\delta_{\bull}))\otimes_{R_{(f)}}R_{(0)} \cong \mathbf{0}.
\end{align*}
From which we deduce that the homology module  $\mathbb{H}^{i}(\mathcal{F}^{\bull}\otimes_{R}R_{(f)},\delta_{\bull})$ is annihilated by a power of $f$, and hence (by \cite{eisenbudbook}  Corollary 2.17 pg. 76) has finite length. \\
\ \\
$ii)$ First observe that \textbf{Det} commutes with localisation:
\begin{align*}
\mathbf{Det}(\mathcal{F}^{\bull}\otimes_{R}R_{(0)},\delta_{\bull})\cong \mathbf{Det}(\mathcal{F}^{\bull}\otimes_{R}R_{(f)},\delta_{\bull})\otimes_{R_{(f)}}R_{(0)}\ .
\end{align*}
Now use \textbf{D2}\footnote{Since $\mathbb{H}^{i}(\mathcal{F}^{\bull}\otimes_{R}R_{(f)},\delta_{\bull})$ is all torsion, we are taking the determinant in the sense of the free resolution $(**)$ . There are many compatibilities to check, but the determinant is insensitive to the choice of  resolution \cite{detdiv}.}
\begin{align*}
\mathbf{Det}(\mathcal{F}^{\bull}\otimes_{R}R_{(f)},\delta_{\bull}) \cong \mathbf{Det}(\mathbb{H}^{i}(\mathcal{F}^{\bull}\otimes_{R}R_{(f)},\delta_{\bull}))\ .
\end{align*}
From $i)$ we know that $\mathbb{H}^{i}(\mathcal{F}^{\bull}\otimes_{R}R_{(f)},\delta_{\bull})$ is a finitely generated \emph{torsion module}, since $R_{(f)}$ is a principal ideal domain, we have, by the structure theorem for such modules, that there are positive integers $n_{i},m_{i1},m_{i2},\dots m_{in_{i}}$ such that
\begin{align*}
\mathbb{H}^{i}(\mathcal{F}^{\bull}\otimes_{R}R_{(f)},\delta_{\bull})\cong \bigoplus _{j=1}^{n_{i}}R_{(f)}/f^{m_{ij}}R_{(f)}.
\end{align*}
From this we get a free resolution $(**)$
\[\begin{CD}
 0@>>>\underset{n_{i}}{\underbrace{\bigoplus}}R_{(f)}@>\delta_{0}>>\underset{n_{i}}{\underbrace{\bigoplus}}R_{(f)}@>\delta_{1}>>\mathbb{H}^{i}(\mathcal{F}^{\bull}\otimes_{R}R_{(f)},\delta_{\bull})@>>> 0 \ . \end{CD} \]
Where the first map is  defined by:
\begin{align*}
 \delta_{0}=
\begin{pmatrix}f^{m_{i1}} & \dots & 0\\
        0&f^{m_{i2}} \dots &  0  \\
         \dots& \dots & \dots \\
          0 & \dots &  f^{m_{in_{i}}}  \\ \end{pmatrix}\ .
\end{align*}

 Therefore $\mathbf{Det}(\mathbb{H}^{i}(\mathcal{F}^{\bull}\otimes_{R}R_{(f)},\delta_{\bull}))$ can be identified with the determinant  of $\delta_{0}$
\begin{align*}
 \mbox{Det}(\delta_{0})=f^{\sum_{1\leq j\leq n_{i}}m_{ij}}\ .
\end{align*}
 Since
 \begin{align*}
 {l_{R_{(f)}}(\mathbb{H}^{i}(\mathcal{F}^{\bull}\otimes_{R}R_{(f)},\delta_{\bull}))}= \sum_{1\leq j\leq n_{i}}m_{ij}\ ,
 \end{align*}
 we may conclude that
\begin{align*}
\mathbf{Det}(\mathbb{H}^{i}(\mathcal{F}^{\bull}\otimes_{R}R_{(f)},\delta_{\bull}))\cong  f^{l_{R_{(f)}}(\mathbb{H}^{i}(\mathcal{F}^{\bull}\otimes_{R}R_{(f)},\delta_{\bull}))}.
\end{align*}
This completes the proof of $ii)$ and so the proposition.

With this in hand we return to the problem of ascertaining the order of vanishing. The main point now is to establish the following
 \begin{theorem}(Cayleys' Theorem on $X$ Resultants) (See \cite{cay},\cite{fischer}, \cite{fog}, \cite{detdiv}, and \cite{ksz})
\begin{align*}\label{cay}
\mathfrak{Tor}(\Gamma(\mathfrak{K}^{\bull}(m)_{X}),\dl_{\bull};\{f_{1,i},f_{2,i}, \dots f_{d_{i},i}\})(l_{0},\dots, l_{n}) = R_{X}((l_{0},\dots, l_{n}))^{(-1)^{n+1}} \ .
\end{align*}
\end{theorem}
In other words, the order of vanishing is just $(-1)^{n+1}$ and hence for $m$ sufficiently large, one has the \emph{canonical} identification of one dimensional complex vector spaces.
\begin{align*}
{\bigotimes_{i=0}^{n+1}{\bf{det}}(H^{0}(X,\mathcal{O}(m+i)))^{(-1)^{i+1}\binom{n+1}{i}}= \mathbb{C}R_{X}^{(-1)^{n+1}}}
\end{align*}

\emph{ Proof}\\
  $\mathfrak{Tor}(\Gamma(\mathfrak{K}^{\bull}(m)_{X}),\dl_{\bull};\{f_{1,i},f_{2,i}, \dots f_{d_{i},i}\})(w_{ij})$ is finite and nowhere vanishing away from $Z= \{R_{X}=0\} $.
Therefore, for all primes $f \neq R_{X}$, $f\in \mathbb{C}[w_{ij}]$ we have
\begin{align*}
ord_{f}(\mathbf{Det}(\mathcal{F}^{\bull}\otimes_{\mathbb{C}[w_{ij}]}\mathbb{C}[w_{ij}]_{(f)}),\delta_{\bull})=0\ , \ \mathcal{F}^{\bull} := H^{0}(X,\mathcal{O}(m+\bull))\otimes \wedge^{\bull}\mathbb{C}^{n+1}\otimes {\mathbb{C}[w_{ij}]} \ .
\end{align*}
 
Therefore we are reduced to computing the order of $R_{X}$ in the determinant. By part $ii)$ of proposition (\ref{order}) this follows immediately from 
\begin{align*}
l_{\mathbb{C}[w_{ij}]_{(R_{X})}}(\mathbb{H}^{i}(\mathcal{F}^{\bull}\otimes_{\mathbb{C}[w_{ij}]}\mathbb{C}[w_{ij}]_{(R_{X})}, \delta_{\bull}))=
\begin{cases}
0&i<n+1\\
1&i=n+1
\end{cases}
\end{align*}
Which completes the proof.  $\Box$
 \begin{center}\textbf{Acknowledgments}\end{center}
The junior author owes an enormous debt to Eckart Viehweg and H\'el\`ene Esnault  for their many invitations to UniversitŠt Duisburg-Essen, his conversations with these two mathematicians were in large part responsible for any clarity that we have achieved in this paper. In particular he would like to thank Eckart for his patient and thorough explanation of Geometric Invariant Theory on the Hilbert Scheme.  
He would also like to thank Paul Gauduchon, Toshiki Mabuchi, Akira Fujiki for very stimulating discussions on this topic in Paris, and Osaka respectively.  
Part of this work was initiated at the Max Planck Institute for Mathematics in the Sciences at Leipzig during the summer of 2005, the first author would like to thank Professor J\"urgen Jost for his kind invitation. 
\bibliography{ref}

\begin{thebibliography}{10}

\bibitem{gauduchon}
Vestislav Apostolov, David~M.J. Calderbank, Paul Gauduchon, and Christina~W.
  Tonnesen-Friedman.
\newblock Hamiltonian 2-forms in kahler geometry, iii extremal metrics and
  stability.
\newblock {\em {math.DG/0511118}}, 2005.

\bibitem{bandmab}
Shigetoshi Bando and Toshiki Mabuchi.
\newblock Uniqueness of {E}instein {K}\"ahler metrics modulo connected group
  actions.
\newblock In {\em Algebraic geometry, Sendai, 1985}, volume~10 of {\em Adv.
  Stud. Pure Math.}, pages 11--40. North-Holland, Amsterdam, 1987.

\bibitem{bredon}
Glen~E. Bredon.
\newblock {\em Sheaf theory}, volume 170 of {\em Graduate Texts in
  Mathematics}.
\newblock Springer-Verlag, New York, second edition, 1997.

\bibitem{cay}
Arthur Cayley.
\newblock On the theory of elimination.
\newblock {\em Cambridge and Dublin Math Journal}, 3, 1848.

\bibitem{geodrays}
X.X. Chen.
\newblock Space of k\"ahler metrics iii--on the lower bound of the calabi
  energy and geodesic distance.
\newblock {\em math.DG/0606228}, 2006.

\bibitem{chentian}
X.X. Chen and G.~Tian.
\newblock Geometry of kahler metrics and foliations by holomorphic discs.
\newblock {\em arXiv}, math.DG/0507148, 2005.

\bibitem{dolga}
Igor Dolgachev.
\newblock {\em Lectures on invariant theory}, volume 296 of {\em London
  Mathematical Society Lecture Note Series}.
\newblock Cambridge University Press, Cambridge, 2003.

\bibitem{skdproj}
S.~K. Donaldson.
\newblock Scalar curvature and projective embeddings. {I}.
\newblock {\em J. Differential Geom.}, 59(3):479--522, 2001.

\bibitem{skdtoric}
S.~K. Donaldson.
\newblock Scalar curvature and stability of toric varieties.
\newblock {\em J. Differential Geom.}, 62(2):289--349, 2002.

\bibitem{eisenbudbook}
David Eisenbud.
\newblock {\em Commutative algebra}, volume 150 of {\em Graduate Texts in
  Mathematics}.
\newblock Springer-Verlag, New York, 1995.
\newblock With a view toward algebraic geometry.

\bibitem{eisenbud}
David Eisenbud, Frank-Olaf Schreyer, and Jerzy Weyman.
\newblock Resultants and {C}how forms via exterior syzygies.
\newblock {\em J. Amer. Math. Soc.}, 16, 2003.

\bibitem{fischer}
E.~Fischer.
\newblock \"{U}ber die {C}ayleysche {E}liminationsmethode.
\newblock {\em Math. Z.}, 26(1):497--550, 1927.

\bibitem{fog}
John Fogarty.
\newblock Truncated {H}ilbert functors.
\newblock {\em J. Reine Angew. Math.}, 234:65--88, 1969.

\bibitem{tableaux}
William Fulton.
\newblock {\em Young tableaux}, volume~35 of {\em London Mathematical Society
  Student Texts}.
\newblock Cambridge University Press, Cambridge, 1997.
\newblock With applications to representation theory and geometry.

\bibitem{harmonicchern}
A~Futaki.
\newblock Harmonic total chern forms and stability.
\newblock {\em math.DG/0603706}, 2006.

\bibitem{gkz}
I.~M. Gelfand, M.~M. Kapranov, and A.~V. Zelevinsky.
\newblock {\em Discriminants, resultants, and multidimensional determinants}.
\newblock Mathematics: Theory \& Applications. Birkh\"auser Boston Inc.,
  Boston, MA, 1994.

\bibitem{globmod}
D.~Gieseker.
\newblock Global moduli for surfaces of general type.
\newblock {\em Invent. Math.}, 43(3):233--282, 1977.

\bibitem{sobolev}
D~Hoffman and J~Spruck.
\newblock Sobolev and isoperimetric inequalities for reimannian submanifolds.
\newblock {\em Comm. Pure Appl. Math.}, 27, 1974.

\bibitem{ksz}
M.~M. Kapranov, B.~Sturmfels, and A.~V. Zelevinsky.
\newblock Chow polytopes and general resultants.
\newblock {\em Duke Math. J.}, 67(1):189--218, 1992.

\bibitem{detdiv}
Finn~Faye Knudsen and David Mumford.
\newblock The projectivity of the moduli space of stable curves {I}:
  {P}reliminaries on ``det'' and ``{D}iv''.
\newblock {\em Math. Scand.}, 39(1), 1976.

\bibitem{mabuchiext}
T.~Mabuchi.
\newblock Stability of extremal k\"ahler manifolds.
\newblock (math.DG/0404211), 2004.

\bibitem{mabuchichow}
T.~Mabuchi.
\newblock The chow-stability and hilbert-stability in mumford's geometric
  invariant theory.
\newblock {\em {math.DG/0607590}}, 2006.

\bibitem{mabuchi}
Toshiki Mabuchi.
\newblock {$K$}-energy maps integrating {F}utaki invariants.
\newblock {\em Tohoku Math. J. (2)}, 38(4):575--593, 1986.

\bibitem{simon}
J.H. Michael and L~Simon.
\newblock Sobolev and mean-value inequalities on generalized submanifolds of
  $\mathbb{R}^{n}$.
\newblock {\em Comm. Pure Appl. Math.}, 26, 1973.

\bibitem{git}
D.~Mumford, J.~Fogarty, and F.~Kirwan.
\newblock {\em Geometric invariant theory}, volume~34 of {\em Ergebnisse der
  Mathematik und ihrer Grenzgebiete (2) [Results in Mathematics and Related
  Areas (2)]}.
\newblock Springer-Verlag, Berlin, third edition, 1994.

\bibitem{sopv}
David Mumford.
\newblock Stability of projective varieties.
\newblock {\em Enseignement Math. (2)}, 23(1-2):39--110, 1977.

\bibitem{ags}
Sean~T. Paul and Gang Tian.
\newblock Analysis of geometric stability.
\newblock {\em Int. Math. Res. Not.}, (48):2555--2591, 2004.

\bibitem{ross&thomas1}
J.~Ross and R.~Thomas.
\newblock An obstruction to the existence of constant scalar curvature k\"ahler
  metrics.
\newblock {\em {math.DG/0412518}}, 2005.

\bibitem{ross&thomas2}
J.~Ross and R.~Thomas.
\newblock A study of the hilbert-mumford criterion for the stability of
  projective varieties.
\newblock {\em {math.AG/0412519}}, 2005.

\bibitem{faisceaux}
Jean-Pierre Serre.
\newblock Faisceaux alg\'ebriques coh\'erents.
\newblock {\em Ann. of Math. (2)}, 61:197--278, 1955.

\bibitem{localalg}
Jean-Pierre Serre.
\newblock {\em Local algebra}.
\newblock Springer Monographs in Mathematics. Springer-Verlag, Berlin, 2000.
\newblock Translated from the French by CheeWhye Chin and revised by the
  author.

\bibitem{psc}
Gang Tian.
\newblock K\"ahler-{E}instein metrics with positive scalar curvature.
\newblock {\em Invent. Math.}, 130(1):1--37, 1997.

\bibitem{eckart}
Eckart Viehweg.
\newblock {\em Quasi-projective moduli for polarized manifolds}, volume~30 of
  {\em Ergebnisse der Mathematik und ihrer Grenzgebiete (3) [Results in
  Mathematics and Related Areas (3)]}.
\newblock Springer-Verlag, Berlin, 1995.

\bibitem{weyman}
Jerzy Weyman.
\newblock {\em Cohomology of vector bundles and syzygies}, volume 149 of {\em
  Cambridge Tracts in Mathematics}.
\newblock Cambridge University Press, Cambridge, 2003.

\end{thebibliography}
\end{document}